\documentclass[11pt]{amsart}
\usepackage{amsmath}
\usepackage{amscd}
\usepackage{amssymb}
\usepackage{pb-diagram}
\usepackage{comment}
\usepackage{graphicx}
\usepackage{color}
\usepackage{enumitem}
\usepackage{float}

\def\refeq#1{\if\workingver y(\ref{#1})-[[#1]]\else(\ref{#1})\fi}
\def\refth#1{\if\workingver y\ref{#1}-[[#1]]\else\ref{#1}\fi}
\def\mylabel#1{\if\workingver y\label{#1}{\bf\ \ [[#1]]\ \ }\else\label{#1}\fi}
\def\mybibitem#1{\if\workingver y\bibitem{#1}{\bf\ \ [[#1]]\ \
}\else\bibitem{#1}\fi}

\newfont{\msam}{msam10}
\newfont{\msbm}{msbm10}

\def\articletheorems{
\newtheorem{thm}{Theorem}[section]
\newtheorem{lem}[thm]{Lemma}

\newtheorem{cor}[thm]{Corollary}
\newtheorem{prop}[thm]{Proposition}

\newtheorem{algo}{Algorithm}[section] 
\newtheorem{alg}[thm]{Algorithm}  

}

\def\oarrow{{\mathrel{\longrightarrow\mkern-25mu\circ}\;\;}}
\def\omap{\oarrow}

\def\cA{\text{$\mathcal A$}}

\def\cK{\text{$\mathcal K$}}
\def\cL{\text{$\mathcal L$}}

\def\cV{\text{$\mathcal V$}}

\newcommand{\id}{\operatorname{id}}
\newcommand{\cl}{\operatorname{cl}}

\newcommand{\bd}{\operatorname{bd}}
\newcommand{\Bd}{\operatorname{Bd}}

\newcommand{\dom}{\operatorname{dom}}

\newcommand{\card}{\operatorname{card}}

\newcommand{\im}{\operatorname{im}}

\renewcommand{\emptyset}{\varnothing}

\def\proof{{\bf Proof:\ }}

\def\begeq#1{\begin{equation}\mylabel{#1}}
\def\endeq{\end{equation}}

\def\mathobj#1{\mbox{$#1$}}

\def\NN{\mathobj{\mathbb{N}}}

\def\RR{\mathobj{\mathbb{R}}}

\def\setof#1{\mbox{$\{\,#1\,\}$}}

\def\0#1{\hbox{\kern25pt}$ #1 $\\}
\def\1#1{\hbox{\kern40pt}$ #1 $\\}
\def\2#1{\hbox{\kern55pt}$ #1 $\\}
\def\3#1{\hbox{\kern70pt}$ #1 $\\}

\newcounter{li}

\def\begalg#1{\begin{algo}\mylabel{#1}\normalshape:\small\baselineskip 10pt\\}
\def\endalg{\end{algo}}

\def\Figures(include=#1,cat=#2){
  \renewcommand{\textfraction}{.20}
  \renewcommand{\topfraction}{.80}
  \renewcommand{\bottomfraction}{.80}
  \renewcommand{\floatpagefraction}{.80}
  \newcount\figcount
  \figcount=0
  \let\includefigures=#1
  \def\figcat{#2}
}

\def\FigureFromFile[#1][#2](#3)#4
{
  \begin{figure}[htbp]
     \global\advance\figcount by 1
     \if\includefigures y\special{anisoscale #1.wmf, \the\hsize #2}\fi
     \vspace{#2}
     \caption{#4}
     \mylabel{#3}
   \end{figure}
}

\def\FigureFromFileTwoD[#1][#2,#3](#4)#5
{
  \begin{figure}[htbp]
     \global\advance\figcount by 1
     \if\includefigures y\special{anisoscale #1.wmf, #2 #3}\fi
     \vspace{#2}
     \caption{#5}
     \mylabel{#4}
   \end{figure}
}

\def\FigureF<#1>[#2](#3)#4
{
  \begin{figure}[htbp]
     \global\advance\figcount by 1
     \if\includefigures y\special{anisoscale \figcat/fig\number\figcount.wmf,
       \the\hsize #2}
     \fi
     \if\includefigures p
       \leavevmode
       \epsfxsize=\hsize
       \epsffile{#1}
     \fi
     \if\includefigures y
          \vspace{#2}
     \fi
     \caption{#4}
     \mylabel{#3}
   \end{figure}
}

\def\Figure[#1](#2)#3
{
  \begin{figure}[htbp]
     \global\advance\figcount by 1
     \if\includefigures y\special{anisoscale \figcat/fig\number\figcount.wmf,
       \the\hsize #1}
     \fi
     \if\includefigures p
       \leavevmode
       \epsfxsize=\hsize
       \epsffile{fig\number\figcount.eps}
     \fi
     \if\includefigures y
          \vspace{#1}
     \fi
     \caption{#3}
     \mylabel{#2}
   \end{figure}
}

\def\0{\hbox{\kern5pt}}
\def\1{\hbox{\kern20pt}}
\def\2{\hbox{\kern35pt}}
\def\3{\hbox{\kern50pt}}
\def\4{\hbox{\kern65pt}}
\def\5{\hbox{\kern80pt}}
\def\6{\hbox{\kern95pt}}

\def\kw#1{\textbf{#1}}

\def\kwif{\kw{if}\;}
\def\kwthen{\;\kw{then}\;}
\def\kwelse{\kw{else}\;}

\def\kwwhile{\kw{while}\;}

\def\kwforeach{\kw{foreach}\;}

\def\kwreturn{\kw{return}\;}
\def\kwdo{\;\kw{do}\;}

\def\kweach{\kw{each}\;}

\def\kwfor{\kw{for}\;}

\def\kweach{\kw{each}\;}

\newcommand{\cbd}{\operatorname{cbd}}

\newcommand{\Cycles}{\operatorname{Cycles}}
\newcommand{\CW}{\operatorname{CW}}
\newcommand{\crit}{\operatorname{crit}}
\renewcommand{\Bd}{\operatorname{Bd}}
\newcommand{\HTW}{\operatorname{HTW}}

\def\Dim#1{|#1|}

\def\opb#1{\stackrel{\circ}{#1}}
\def\clb#1{#1}
\def\cell#1{\langle #1 \rangle}

\def\mathobj#1{\mbox{$#1$}}
\def\RR{\mathobj{\mathbb{R}}}
\def\c#1{\mathobj{\mathcal{#1}}}

\def\setof#1{\mbox{$\{\,#1\,\}$}}

\articletheorems

\begin{document}

    \author[P. Brendel]{Piotr Brendel}
    \address{Piotr Brendel, Division of Computational Mathematics,
      Jagiellonian University in Krakow, Poland.
      E-mail: Piotr.Brendel@ii.uj.edu.pl
    }
    \author[G. Ellis]{Graham Ellis}
    \address{Graham Ellis,
      School of Mathematics, National University of Ireland, Galway, Ireland.
      E-mail: graham.ellis@nuigalway.ie
    }
    \author[M. Juda]{Mateusz Juda}
    \address{Mateusz Juda, Division of Computational Mathematics,
      Jagiellonian University in Krakow, Poland.
      E-mail: Mateusz.Juda@ii.uj.edu.pl
    }
    \author[M. Mrozek]{Marian Mrozek}
    \address{Marian Mrozek, Division of Computational Mathematics,
      Jagiellonian University in Krakow, Poland.
      E-mail: Marian.Mrozek@ii.uj.edu.pl
    }

\title[Fundamental Group Algorithm]{Fundamental Group Algorithm for low dimensional tessellated CW complexes.}
\thanks{This research is partially supported
  by the EU under the {\sc Toposys} project FP7-ICT-318493-STREP,
  and by the ESF under the ACAT Research Network Programme.}

\begin{abstract}
We present a detailed description
of a fundamental group algorithm based on Forman's combinatorial version of Morse theory. We use this algorithm in a classification problem of prime knots up to $14$ crossings.
\end{abstract}

\maketitle

\section{Introduction}

The aim of this paper is twofold. We present the mathematical details of the algorithm for computing a presentation of the fundamental group of a finite regular CW complex which was introduced in~\cite{Brendel-etal-2014} and which is based on Forman's combinatorial version of Morse theory~\cite{Fo98a}. We also present details of the performance of a C++ implementation of the algorithm. To assess the performance we systematically apply the algorithm to spaces arising as complements of prime knots on fourteen or fewer crossings and use the resulting fundamental group presentation, with standard low-index subgroup procedures, to distinguish between the knots.

\subsection{Prior work}
The combinatorially defined edge-path group of a connected simplicial complex K, due to
Reidemeister, is well-known to be isomorphic to the fundamental group $\pi_1(K)$ (see \cite{Spanier-1981}).
It is also
well-known that this combinatorial definition and isomorphism extends to connected regular
CW-spaces. In this paper we use the terms edge-path group and fundamental group interchangeably as
synonyms. Several authors have described algorithms for implementing Reidemeister's edge-path
group. Rees and Soicher \cite{ReesSoicher-2000} use spanning trees and redundant relator searches in their description
of an algorithm for finding a small finite presentation of the edge-path group of a 2-dimensional
combinatorial cell complex; they implement their algorithm in {\sc GAP} \cite{gap} for 2-dimensional simplicial
clique complexes of graphs. Letscher \cite{Letscher-2012} uses spanning trees and Tietze elimination/reduction of
relators to compute edge-path groups from the 2-skeleta of simplicial complexes arising from knot
complements, the knots being produced from experimental data on protein backbones. Palmieri
et al. \cite{Palmieri-2009} have implemented the edge-path group of simplicial complexes in Sage \cite{Stein-2013};
the implementation uses the 2-skeleton of the complex and calls {\sc gap}'s Tietze reduction/elimination
procedures \cite{gap}. Kim et al. \cite{Kim-etal-2008} describe an algorithm for the fundamental groups of 3-dimensional
simplicial complexes; their algorithm, which makes use of 3-dimensional cells and the language of
general CW-spaces, is applied to 3-dimensional tetrahedral meshes arising in computer vision.

We also mention that there is a large literature on the related problem of algorithmically
determining a collection of shortest generating loops for the fundamental group of a space. The
case of oriented combinatorial 2-manifolds is treated by Ericson and Whittlesey in [10]; their
algorithm involves the computation of spanning trees.

\subsection{Our contribution}
Recall that a combinatorial vector field on a regular CW complex $X$
is a partition $\cV$ of cells of $X$ into singletons and doubletons such that
each doubleton consists of a
cell and one of its facets, that is, one of its faces of codimension 1.
The singletons are called {\em critical cells}, the doubletons are called {\em vectors}.
The combinatorial vector field is {\em acyclic} if the facet digraph of X, with direction reversed on edges in $\cV$, is acyclic.
The algorithm we propose (Algorithm~\ref{alg:discreteMorseFunction})  inputs a connected finite regular CW complex $X$,
constructs a homotopy equivalent CW complex $X'$ with potentially fewer cells,
and outputs the finite presentation for $\pi_1(X)$ corresponding to the 2-skeleton of $X'$.
The algorithm runs in linear time with respect to the number of cells in $X$ (see Theorem~\ref{thm:discreteVectorField-algo}).
However, its usefulness depends on the size of $X'$. We have no estimate for the number of cells in $X'$ other than it will be no greater than the number of cells in $X$.
However, from several numerical experiments we see that  $X'$ is usually significantly smaller than $X$.

At the heart of the algorithm is Forman's theorem \cite[Corollary 3.5]{Fo98a}
stating that an acyclic combinatorial
vector field $\cV$ on a regular CW complex $X$ leads to another,
homotopic, CW complex whose cells are in one-to-one correspondence
with the critical cells of $\cV$.
In order to justify the correctness of the algorithm
we rework in an algorithmic, recursive spirit
the proof of Forman's theorem. We actually show that the vectors
of $\cV$ may be quotiented out one by one until no vector is left.

We implemented the algorithm in C++ \cite{RedHom}.
The implementation aims at large CW complexes arising from real-world data such as protein backbone
complements or point-cloud data sets.
However, in order to systematically study the performance of the algorithm, in this paper we apply it to a family of small cubical complexes arising as the complements of prime knots given in terms of efficient arc presentations \cite{KnotAtlas}. We have chosen to work with this case,
because it provides a large source of examples whose complexity can, to some extent, be measured by the number of crossings in the knots. Recall that in  \cite{Brendel-etal-2014} we proposed a computable,
algebraic invariant of prime knots and used it to classify,
up to mirror image, all prime knots up to 11 crossings.
The invariant is constructed from low-index subgroups of the fundamental group of the complement
of the knot with the index not exceeding $6$.
The C++ implementation of our algorithm enables fast computation of the presentation
of the fundamental group. This lets us increase the classification in \cite{Brendel-etal-2014}
from  $801$ knots up to $11$ crossings to $59937$ prime knots up to $14$ crossings
in the tabulation provided by Hoste, Thistlethwaite and Weeks in \cite{Hoste-etal-1998}
and available from \cite{KnotAtlas}.
The C++ implementation shifted the computational barrier from the presentation of the fundamental group computations to
the algebraic invariant computations. The latter are performed in GAP \cite{gap} and HAP \cite{Ellis-2013}.

The organization of the paper is as follows.
In Section~\ref{sec:preliminaries} we gather notation and preliminaries.
In Section~\ref{sec:disc-morse} we present a recursive approach to the discrete Morse theory.
In Section~\ref{sec:reductions} we study various reduction techniques.
In Section~\ref{sec:red-algo} we present the algorithms.
In Section~\ref{sec:applications} we apply the results of the paper to knot classification.

\section{Preliminaries}
\label{sec:preliminaries}

\subsection{Notation}

Let $\NN:=\{1,2,3,\ldots\}$ stand for the set of natural numbers, $\RR$
for the set of reals and let $I$ denote the interval $[0,1]$.
We write $f:X\omap Y$ for a {\em partial map} from $X$ to $Y$, that is a map
defined on a subset $\dom f\subset X$, called the {\em domain} of $f$,
and such that the set of values of $f$, denoted $\im f$, is contained in $Y$.
Given a topological space $X$ we write $\cl A$ for the closure of $A\subset X$.

For $n\in\NN$ we write
\begin{eqnarray*}
    \clb{B^n}&:=&\setof{x\in\RR^n\mid\;\;||x||\leq 1},\\
    S^{n-1}&:=&\setof{x\in\RR^n\mid\;\;||x||=1},\\
    B^{n-1}_-&:=&\setof{x=(x_1,x_2,\ldots x_n)\in S^{n-1}\mid\;\;x_n\leq 0},\\
    B^{n-1}_+&:=&\setof{x=(x_1,x_2,\ldots x_n)\in S^{n-1}\mid\;\;x_n\geq 0},\\
    \opb{B^n}&:=&\setof{x\in\RR^n\mid\;\;||x||<1}.
\end{eqnarray*}
Note that in particular $\clb{B^1}=[-1,1]$ and $S^0=\{-1,1\}$.
We extend this notation to $n=0$ by setting
\begin{eqnarray*}
    \opb{B^0}&:=&\clb{B^0}:=\{\ast\},\\
    S^{-1}&:=&\emptyset,
\end{eqnarray*}
where $\ast$ is a unique point not belonging to any other set considered in this paper.

Let $A$ be any set. An {\em oriented element} of $A$ is $a^\epsilon :=(a,\epsilon)$
for $\epsilon\in\{-1,1\}$ and $a \in A$. We denote the set of oriented elements of $A$ by $\hat{A}$.
We identify $a^1$ with $a$ and we say that $a$ and $a^{-1}$ have mutually inverse orientation.

By a {\em path} in a topological space $X$ we mean a continuous function $\theta:\clb{B^1}=[-1,1]\to X$.
We identify the oriented path $\theta^1$ with $\theta$ and $\theta^{-1}$ with the path
\[
  \clb{B^1}\ni t\mapsto \theta(-t)\in X.
\]
We refer to $\theta^{-1}$ as the {\em inverse path} of $\theta$.

\subsection{Directed multigraphs}

Recall that a {\em directed multigraph} is a quadruple $G=(V,E,e_{-1},e_{+1})$, where
$e_{-1},e_{+1}:E\to V$ are given maps. The set $V$ is the set of {\em vertices} of $G$,
the set $E$ is the set of {\em edges} of $G$ and the maps $e_{-1},e_{+1}$ assign
to each edge its {\em initial} and {\em terminal} vertex respectively.

A {\em path} in $G$ is a sequence of oriented edges $\eta:=e^{\epsilon_1}_1e^{\epsilon_2}_2\cdots e^{\epsilon_n}_n$
such that $e_{\epsilon_k}(e_k)=e_{-\epsilon_{k+1}}(e_{k+1})$ for $k=1,\ldots n-1$.
The {\em initial vertex} of the path $\eta$ is $e_{-\epsilon_{1}}(e_{1})$ and the {\em terminal vertex} of the path $\eta$ is $e_{\epsilon_{n}}(e_{n})$.
The path $\eta$ is a {\em cycle} if the terminal and initial vertices coincide.
The {\em inverse of the path} $\eta$ is $\eta^{-1}:=e^{-\epsilon_n}_ne^{-\epsilon_{n-1}}_{n-1}\cdots e^{-\epsilon_1}_1$.

Let $\zeta$ be a path in $G$ with the same initial and terminal vertices as an edge $z\in E$.
Moreover, assume that neither $z$ nor $z^{-1}$ appear in $\zeta$. Let $\eta$ be a path in $G$.
A {\em $(z,\zeta)$-substitution in $\eta$}, denoted $S_{z,\zeta}(\eta)$, is the path
obtained from $\eta$ by replacing every occurrence of $z$ in $\eta$ by $\zeta$ and every occurrence
of $z^{-1}$ in $\eta$ by $\zeta^{-1}$.

\subsection{Quotient spaces}

We recall from \cite{Engelking-1989} some facts concerning quotient topological
spaces.
Given an equivalence relation $R$ on a topological space $X$,
we denote by $[x]_R$ the equivalence class of $x\in X$ and by
$X/R$ the set of all equivalence classes of $R$. Let
$\kappa_R:X\ni x\mapsto [x]_R\in X/R$ be the canonical projection.
We drop the subscript $R$ whenever the relation $R$ is clear from context.
The {\em quotient topology} on $X/R$ consists of
subsets $A\subset X/R$ such that $\kappa^{-1}(A)$ is open in $X$.

\begin{prop}
\label{prop:quotient-iterated}
\cite[Proposition 2.4.14]{Engelking-1989}
Assume $R$ is an equivalence relation on $X$ and $S$ is an equivalence relation on $X/R$.
Then, $X/R/S$ is homeomorphic to $X/R_S$, where $R_S$ is an equivalence relation on $X$
defined by $xR_Sy$ if and only if $\kappa_S(\kappa_R(x))=\kappa_S(\kappa_R(y))$.
\end{prop}

Given a partial continuous map $f:Y\oarrow X$ such that $\dom f$ is closed in $Y$
we define $X\cup_f Y$ as the quotient space
$(X\sqcup Y)/R_f$ of the disjoint union $X\sqcup Y$ by the equivalence relation $R_f$
which identifies $x\in\dom f$ with $f(x)\in X$. For $z\in X\sqcup Y$ we denote by $[z]_f$
the equivalence class of the relation $R_f$.
Note that if $f$ is empty, then $X\cup_f Y$ coincides with $X\sqcup Y$.
In the special case when $Y=B^n$ and $\dom f=S^{n-1}$ we say that
$X\cup_f Y$ results from $X$ by {\em gluing  an $n$ dimensional cell} via the {\em attaching map} $f$.

The relation $R$ is called {\em closed} if the map $\kappa_R$ is closed.
The following proposition is a special case of \cite[Proposition 2.4.9]{Engelking-1989}.
\begin{prop}
\label{prop:quotient-usc}
The relation $R$ is closed if and only if for any open $U$ in $X$ the set
\[
  U_R:=\setof{x\in U\mid [x]_R\subset U}
\]
is open in $X$.\qed
\end{prop}

The following proposition is an easy consequence of Proposition~\ref{prop:quotient-usc}
\begin{prop}
\label{prop:quotient-extension-usc}
Assume $Y\subset X$ is closed and $R$ is a closed equivalence relation in $Y$.
Then $\bar{R}:=R\cup \id_X$ is a closed equivalence relation on $X$.\qed
\end{prop}

\begin{prop}
\label{prop:quotient-Hausdorff}
Assume $X$ is a Hausdorff topological space and $R$ is a closed equivalence relation
in $X$ with compact equivalence classes. Then $X/R$ is a Hausdorff space.
\end{prop}
\proof
  Let $[x]\neq [y]$ be two equivalence classes of $R$. Since they are compact and $X$
  is a Hausdorff space, we can find disjoint open sets $U$, $V$ in $X$ such that $[x]\subset U$
  and $[y]\subset V$. By Proposition~\ref{prop:quotient-usc} $U_R$ and $V_R$ are open in $X$.
  It follows that $\kappa(U_R)$, $\kappa(V_R)$ are disjoint, open in $X/R$.
  Obviously, they are neighbourhoods respectively of $[x]$ and  $[y]$.
  Hence, $X/R$ is a Hausdorff space.
\qed

Given $A\subset X$ we denote by $X/A$ the quotient space $X/R_A$, where
\[
  R_A:=\setof{(x,y)\in X^2\mid \text{$x=y$ or $x\in A$ and $y\in A$.}}
\]

\subsection{Finite CW complexes}
 In this paper by a CW complex we always mean a finite CW complex. Here we recall its definition.
Let $K$ be a Hausdorff topological space. A {\em finite CW structure on $K$} is a pair
$(\cK,\{\varphi_\sigma\}_{\sigma\in\cK})$ such that $\cK$ is a finite family of subsets of $K$
and:
\begin{itemize}
   \item[(CW0)] Each element $\sigma\in\cK$ is a subset of $K$ homeomorphic to $\opb{B^n}$ for some $n\geq 0$.
   The subset $\sigma$ is referred to as an {\em $n$-cell}. The number $n$ is called the {\em dimension} of $\sigma$
   and denoted $\Dim{\sigma}$. The dimensions provide a {\em filtration} $K^j:=\bigcup\setof{\sigma\in\cK\mid \Dim{\sigma}\leq j}$.
   \item[(CW1)] The family $\cK$ is a decomposition of $K$, i.e. $K=\bigcup\cK$ and any two different elements of $\cK$ are disjoint.
   \item[(CW2)] Each $\varphi_\sigma$  is a continuous map
   $\varphi_\sigma: \clb{B^{\Dim{\sigma}}} \to K$ mapping
   $\opb{B^{\Dim{\sigma}}}$ homeomorphically  onto $\sigma$ and $\varphi_\sigma(S^{{\Dim{\sigma}}-1})\subset K^{{\Dim{\sigma}}-1}$.
   The map $\varphi_\sigma$ is called the {\em characteristic map} of $\sigma$.
\end{itemize}
A {\em finite CW complex} (cf. \cite[Section I.3]{Cohen-1973}) is a triple $(K,\cK,\{\varphi_\sigma\}_{\sigma\in\cK})$ such that $(\cK,\{\varphi_\sigma\}_{\sigma\in\cK})$ is a finite CW structure on $K$.
It is convenient to slightly abuse the language and refer to $K$ as the CW complex, assuming that the CW structure is given implicitly.

The CW complex is called {\em regular} if each $\varphi_\sigma: \clb{B^{\Dim{\sigma}}} \to K$
is a homeomorphism onto its image.

By the {\em boundary} of a cell $\sigma$ we mean the set $\Bd\sigma:=\cl\sigma\setminus\sigma$. We say that a cell $\tau$ is a {\em face} of a cell $\sigma$ if $\tau\cap\cl\sigma\neq\emptyset$. The face is {\em proper} if $\tau\neq\sigma$.
The face is a {\em facet} if $\Dim{\tau}=\Dim{\sigma}-1$.
A $0$-dimensional cell is called a {\em vertex}.

We refer to the collection of facets of a cell $\sigma$ as its {\em combinatorial boundary}
and we denote it by $\bd\sigma$. It is also convenient to consider the collection
of cells whose facet is $\sigma$. This collection is called
the {\em combinatorial coboundary} of $\sigma$ and is denoted $\cbd\sigma$.
Note that if $\sigma$ is a $1$-cell,
then $\bd\sigma=\{\varphi(-1), \varphi(1)\}$.
In the case of a $1$-cell these two facets (or one if $\varphi(-1)=\varphi(1)$) are referred to as the {\em endpoints} of $\sigma$.

A cell is {\em top-dimensional} if it is not a proper face of a higher dimensional cell.
The CW structure and CW complex are called {\em pure} if all top-dimensional cells have the same
dimension. We will refer to the closure of a top-dimensional cell as a {\em toplex}.
It is easily seen that a CW complex is the union of its toplexes.

A union $L:=\bigcup\cL$ of a subfamily $\cL\subset \cK$ such that $L$ is closed  in $K$
is called a {\em subcomplex} of $K$. A subcomplex is a CW complex with its CW structure consisting
of  $\cL$ with the same characteristic maps as in $\cK$.
Note that, in particular, every vertex is a subcomplex.
Also, $\bigcup\cK^j$ is a subcomplex, called the {\em $j$-skeleton} of $K$.

Since by \cite[Theorem 1.4.10]{FritschPiccinini-1990} the closure of a cell
in a regular CW complex is a subcomplex, we have the following
proposition.
\begin{prop}
\label{prop:intersection-implies-face}
If cells $\sigma,\tau$ in a regular CW complex satisfy
$\tau\,\cap\,\cl\sigma\neq\emptyset$, then $\tau\subset\cl\sigma$,
i.e. $\tau$ is a face of $\sigma$.\qed
\end{prop}

As a consequence of Proposition~\ref{prop:intersection-implies-face}
we obtain the following proposition.
\begin{prop}
\label{prop:facets}
Assume $\tau$ is a proper face of $\sigma$ in a regular CW complex $X$.
Then $\tau$ is a facet of a face of $\sigma$.
\end{prop}
\proof
Set  $n:=\Dim{\sigma}$. Since $X$ is regular, $\Bd\sigma$
is homeomorphic to $S^{n-1}$.
Let $k:=\Dim{\sigma}-\Dim{\tau}$. Since $\tau$ is a proper face of $\sigma$
we have $k>0$. We proceed by induction on $k$. If $k=1$ the claim is obvious.
Thus, assume $k>1$. Let $x\in\tau$ and consider a shrinking sequence
of open balls in $\Bd\sigma$ centered at $x$.
All these balls are homeomorphic to $\opb{B^{n-1}}$ and none of these balls is
contained in $\tau$, because otherwise $\Dim{\tau}=n-1$ and $k=1$.
Thus, there exists a sequence $\{x_j\}\subset\Bd\sigma\setminus\tau$
converging to $x$. By passing to  a subsequence, if necessary,
we may assume that the sequence is contained in a cell $\sigma'\subset\Bd\sigma$.
It follows that $\tau\cap\cl\sigma'\neq\emptyset$ and, by the regularity of $X$,
$\tau$ is a face of $\sigma'$. The conclusion now follows from the induction
assumption, because $\Dim{\sigma'}-\Dim{\tau}<k$.
\qed

\begin{prop}
\label{prop:cw-metrizable}
(\cite[Proposition 1.5.17]{FritschPiccinini-1990})
A finite CW complex is a metrizable topological space.\qed
\end{prop}

Note that if $h:\clb{B^n}\to \clb{B^n}$ is a homeomorphism, then $h$ maps homeomorphically
$\opb{B^n}$ onto $\opb{B^n}$ and $S^{n-1}$ onto $S^{n-1}$. In consequence, if $\varphi_\sigma$ is a characteristic
map, then so is $\varphi_\sigma \circ h$. Thus, we can freely substitute homeomorphisms into
characteristic maps without loosing the CW structure. In particular, in the characteristic map
$\varphi_\sigma$ we can replace $\clb{B^n}$ by any set homeomorphic to $\clb{B^n}$.

Let $\sigma$ be a top-dimensional cell of $K$. Then $K\setminus\sigma$
is easily seen to be a subcomplex of $K$ and $K$ results from
$K\setminus\sigma$ by gluing $B^{\Dim{\sigma}}$
via the attaching map
$\theta_\sigma:=\varphi_{\sigma|S^{\Dim{\sigma}-1}}$.

Let $\sigma=(\sigma_j)_{j=1}^k$ be an {\em admissible ordering} of  the cells of $\cK$,
i.e. an ordering of all cells in $\cK$ such that each cell is preceded by its faces.
For $j=1,2,\ldots k$ set
\[
  \Sigma_j:=\bigcup_{i=1}^j\sigma_i.
\]
It is not difficult to prove  that each $\Sigma_{j}$ is closed in $K$, hence a subcomplex of $K$ and
$\Sigma_{j}$ results from $\Sigma_{j-1}$ by gluing the cell $\sigma_j$ via the attaching map
$\theta_{\sigma_j}=\varphi_{\sigma_j|S^{\Dim{\sigma_j}-1}}$.

One easily verifies that an admissible ordering always exists on a finite CW complex.
Moreover, when $L$ is a subcomplex of $K$, then one can choose the ordering in such a way
that $L=\Sigma_l$ for some $l\leq k$. Thus, we have the following standard result.

\begin{prop}
\label{prop:glueing}
Each finite CW-complex can be constructed recursively by attaching a cell to a subcomplex. \qed
\end{prop}

We will use also the following fact which is a special case of \cite[Theorem 11.11]{Kozlov-2008}.
\begin{thm}
\label{thm:kozlov}
Assume $X_1$, $X_2$ are homotopy equivalent spaces and $h:X_1\to X_2$ is a homotopy equivalence.
If maps $f_i:S^{n-1}\to X_i$ are such that $f_2=h\circ f_1$, then the spaces $X_1\cup_{f_1}B^n$
and $X_2\cup_{f_2}B^n$ are homotopy equivalent.\qed
\end{thm}

\subsection{Oriented CW complexes}
Our  fundamental group algorithm of a CW complex $X$
is based on a theorem reducing the computation to a group generated
by the $1$-cells with relators provided by the $2$-cells. To present the theorem in detail,
we recall the definition of the {\em orientation} of an $n$-cell $\sigma$ in $X$.
For $n=0$ it is just an element of $\{-1,1\}$.
For $n>0$ it is an equivalence class in the homotopy equivalence relation of
the characteristic map $\varphi_\sigma:(\clb{B^n},S^{n-1})\to(\cl\sigma,\Bd\sigma)$
considered as a map of pairs.
By \cite[Corollary 2.5.2]{Geoghegan-2008}
there are precisely two orientations for each cell.
An {\em oriented CW complex} is a CW complex with a fixed orientation for
every its cell.

Assume $K$ is an oriented $2$-dimensional CW complex.
By an {\em edge} in $K$ we mean an oriented $1$-cell.
An edge $\tau$, as each $1$-cell, has two endpoints $\varphi_\tau(-1)$ and $\varphi_\tau(1)$. The orientation
allows us to distinguish between them.
The first is called the initial vertex
of $\tau$ and denoted $e_-(\tau)$,
the other is called the terminal vertex of $\tau$ and denoted $e_+(\tau)$. Thus, the CW structure of the $1$-skeleton
of $K$ may be considered as a directed multigraph.
A characteristic map of an edge is a characteristic
map $\varphi_\tau$ of the associated $1$-cell
which induces the chosen orientation of the edge.

Consider now a $2$-cell $\sigma$. There is a closed path
$\tau_1,\tau_2,\ldots,\tau_n$ in the directed multigraph of edges
such that the homotopy class of the attaching map $\theta_\sigma$
contains a map $\theta:S^1\to K^1$ such that
$S^1=\bigcup_{j=1}^nI_j$
and $\theta_{|I_j}$ is the characteristic map of the edge $\tau_j$
(see \cite[Proposition 3.1.1 and 3.1.2]{Geoghegan-2008}).
The loop $\tau_1,\tau_2,\ldots,\tau_n$ is called the homotopical
boundary of $\sigma$.
We have the following theorem.
\begin{thm}
\label{thm:edge-group} (\cite[Proposition 3.1.7 and Theorem 3.1.8]{Geoghegan-2008})
Assume $K$ is a connected CW complex with precisely one vertex.
Then, the fundamental group of $K$ depends only
on the $2$-skeleton of $K$. Moreover, up to an isomorphism it
is the group generated by the edges of $K$ with arbitrarily selected
orientation of $K$ and homotopical boundaries of all $2$-cells as relators.
\end{thm}

\section{Discrete Morse theory}
\label{sec:disc-morse}

\subsection{Reduction pairs.}
We say that a facet $\alpha_0$ of an $n$-cell  $\alpha_1$ is {\em regular} if
$\varphi_{\alpha_1}$, the characteristic map of $\alpha_1$, satisfies
\begin{itemize}
   \item[(i)]   $\varphi_{\alpha_1}(S^{n-1})$ is a subcomplex of $K$,
   \item[(ii)]  $\varphi_{\alpha_1|\clb{B^n}_-}$ is the characteristic map of $\alpha_0$,
   \item[(iii)] $\varphi_{\alpha_1}(\clb{B^n}_+)\cap\alpha_0=\emptyset$.
\end{itemize}

Note that in a regular CW complex each facet is automatically regular.

We say that $\alpha:=(\alpha_0,\alpha_1)$ is a {\em reduction pair} in a CW complex $K$
if $\alpha_0$ is a regular facet of $\alpha_1$.
The following proposition is straightforward.

\begin{prop}
\label{prop:sub-reduction-pair}
If $L$ is a subcomplex of $K$ and $\alpha$ is a reduction pair in $L$, then $\alpha$ is a reduction pair in $K$. \qed
\end{prop}

\subsection{Discrete vector fields.}
A {\em discrete vector field} $\cV$ on the CW structure $\cK$
is a partition of $\cK$ into singletons and doubletons such that each doubleton,
when ordered by dimension, forms a reduction pair.
By a {\em critical cell} of $\cV$ we mean a cell $\alpha\in\cK$ such that $\{\alpha\}$
is a singleton in $\cV$. We denote the set of critical cells of $\cV$ by $\crit\cV$.
By pairing the lower dimensional cell with the higher dimensional cell in every doubleton of $\cV$
we may and will consider $\cV$ as an injective partial self-map $\cV:\cK\omap\cK$ such that
$\cK$ is the disjoint union of $\dom \cV$, $\im \cV$ and $\crit\cV$
and for each $\alpha\in\dom\cV$ the pair $(\alpha,\cV(\alpha))$ is a reduction pair.
We will call such a reduction pair a {\em vector} of $\cV$.

With $\cV$ we associate a directed graph, $G_\cV$. Its vertices are the cells of $\cK$ and edges go from
each cell $\sigma$ to each of its facets $\tau$ with orientation reversed if $\sigma=\cV(\tau)$.
If $G_\cV$ is acyclic, the discrete vector field $\cV$ is called {\em acyclic}. In this case the graph $G_\cV$ induces a partial order
on the cells in $\cK$. In particular, for each cell $\sigma\in\cK$
there is a minimal cell less than or equal to $\sigma$ in this partial order.
We say about such a cell that it is {\em inferior} to $\sigma$.

Assume now that $K$ is a regular CW complex with a CW structure $\cK$.
An algorithm constructing an acyclic discrete vector field on $\cK$
is presented in Table~\ref{tab:discreteVector Field}.
It is a simplified version of the algorithm proposed
in \cite{Harker-etal-2010,Harker-etal-2014}, based
on the method of coreductions \cite{MrozekBatko-2008}.

\begin{table}
\begin{tabular}{|c|}
\hline
~~\\
\begin{minipage}{120mm}
\begin{algo}\upshape
\label{alg:discreteMorseFunction}
discreteVectorField(\texttt{CW structure} \cK)\\
\0{$V:=L:=C:=\emptyset$;} \\
\0{\kwwhile $\cK\neq\emptyset$ \kwdo} \\
\1{  \kwif $L=\emptyset$\kwthen} \\
\2{    $n:=\min\setof{q\mid \cK^q\neq\emptyset}$;}\\
\2{    move an element $r$ from $\cK^n$ to $C$;}\\
\2{    \kwfor\kweach $u\in(\cK\cap\cbd r)\setminus L$ \kwdo enqueue $u$ in $L$;}\\
\1{  \kwelse} \\
\2{    $\sigma$:=dequeue$(L)$;}\\
\2{    \kwif $\cK\cap\bd \sigma=\{\tau\}$  \kwthen}\\
\3{      remove $\tau$ and $\sigma$ from $\cK$ and insert $(\tau,\sigma)$ to $V$;}\\
\3{      \kwfor\kweach $u\in(\cK\cap\cbd \tau)\setminus L$ \kwdo enqueue $u$ in $L$;}\\
\2{    \kwelse \kwif $\cK\cap\bd \sigma=\emptyset$ \kwthen}\\
\3{      \kwfor\kweach $u\in(\cK\cap\cbd \sigma)\setminus L$ \kwdo enqueue $u$ in $L$;}\\
\0{\kwreturn V;}
\end{algo}
\end{minipage}\\
~~\\
\hline
\end{tabular}
\vspace*{3mm}
\caption{Discrete Vector Field by coreduction method}
\label{tab:discreteVector Field}
\end{table}

For a CW structure $\cK$ let $M(\cK)$ denote the maximal
cardinality of $\bd\sigma$ and $\cbd\sigma$ over all $\sigma\in\cK$.

\begin{thm}
\label{thm:discreteVectorField-algo}
Assume Algorithm~\ref{alg:discreteMorseFunction} in Table~\ref{tab:discreteVector Field}
is called with $\cK$ containing a CW structure
of a non-empty, regular, connected finite CW complex $K$.
The algorithm always stops, returning
an acyclic vector field on $\cK$.
Moreover, on return $C$ contains the respective set of critical cells
with precisely one critical cell in dimension zero.
If  the data structure used to store the CW complex provides
iteration over the boundary and coboundary of a cell as well as removing a cell in $O(M)$ time
with $M:=M(\cK)$, then the algorithm runs in  $O(M^2\card\cK)$ time.
\end{thm}
\proof
Observe that if an element $\tau\in \cK$ contributes its coboundary elements to
$L$ as an element of a coreduction pair $(\tau,\sigma)$, then it is removed from
$\cK$, so it may never again contribute in this role. Similarly,
if an element $\sigma$ with zero boundary contributes its coboundary elements
to $L$, it may never again do so, because it is not in the
coboundary of any element, so it may not reappear in $L$. Therefore, each
element of $S$ may appear in $L$ at most $2M$ times.
It follows that the $\kwwhile$ loop is passed at most
$2M\card \cK$ times.
In particular, the algorithm always stops and its complexity is
$O(M^2\card\cK)$.

Let $\cK_j$, $C_j$ and $V_j$ denote respectively
the contents of variable $\cK$, $C$ and $V$ on leaving the $j$th pass of the $\kwwhile$  loop for $j>0$ and the initial values of these variables for $j=0$. We will show by induction in $j$ that
$V_j$ is an acyclic vector field on $\cK$ with $\cK_j\cup C_j$ as
the set of its critical cells.
This is obvious for $j=0$, because $V_0=\emptyset$, so every cell is critical. Thus, assume the claim holds for $k<j$.
It is obvious that $V_j$ is a discrete vector field on $\cK$.
We only need to prove that $V_j$ is acyclic.
This is evident if $V_j=V_{j-1}$, thus assume $V_j\neq V_{j-1}$.
Then $V_j=V_{j-1}\cup\{(\tau,\sigma)\}$ for some reduction pair $(\tau,\sigma)$.
Assume $V_j$ is not acyclic. Then, there exists a cycle
$s_0,t_0,s_1,t_1,\ldots s_n,t_n,s_0$ in $G_{V_j}$
such that $(t_i,s_i)\in V_j$ and $t_i$ is a facet of $s_{(i+1)\bmod n}$.
In particular, $t_n\in\bd s_0$.
Since by induction assumption $V_{j-1}$ is acyclic,
precisely one of the pairs $(t_i,s_i)$ is the pair appearing
in $V_j$ on the $j$th pass of the $\kwwhile$ loop.
Without loss of generality we may assume that
this is the pair $(t_0,s_0)$.
Thus, we have $\cK_{j-1}\cap\bd s_0=\{t_0\}$.
In particular, $t_0\in\cK_{j-1}$.
Since the pair $(t_1,s_1)\in V_{j-1}$, it must have been included
in the vector field on an earlier  pass of the \kwwhile loop,
say  $k$th pass with $k<j$. Thus, $\cK_{k-1}\cap\bd s_1=\{t_1\}$.
But $t_0\in\bd s_1$, therefore $t_0\not\in\cK_{k-1}$.
However, the sequence $(\cK_i)$ is clearly decreasing, hence
$t_0\not\in\cK_{j-1}$, a contradiction. This proves that $V_j$
is acyclic for every $j$. In particular, the vector field returned
by the algorithm is acyclic. Since after the last pass of the
$\kwwhile$ loop $\cK=\emptyset$, its critical cells are all stored
in $C$.

We still need to prove that $C$ contains precisely one vertex.
First observe that for each
$j$ the union of the cells in  set $\cA_j:=C_j\cup\dom V_j\cup\im V_j$ is a closed subcomplex of $K$.
Indeed, this is obviously true for $j=0$, because $\cA_0=\emptyset$.
Now, arguing by induction, we see that $\cA_{j+1}\setminus \cA_j$
is either empty, or a critical cell in $C_{j+1}$
whose boundary is in $\cA_j$, or
a doubleton $\{\tau,\sigma\}$ such that $(\tau,\sigma)$ is a reduction pair
satisfying $\cK_{j}\cap\bd \sigma=\{\tau\}$, i.e. $\bd \sigma\setminus \cA_{j}=\{\tau\}$. Thus, in each case $\bigcup \cA_j$ is closed,
because $\bigcup \cA_{j-1}$ is closed by induction argument.

Observe that on the very first
pass of the $\kwwhile$ loop a vertex, say $v$,
is removed form $\cK$ and moved
to $C$. Thus, we only need to show that this never happens again.
Assume the contrary. Let $w$ be another vertex moved to $C$
and assume this happens on the $j$th pass of the $\kwwhile$ loop.
Since $K$ is connected, there is a path
$v=v_0,e_0,v_1,e_1,\ldots e_m,v_m=w$ such that each $e_i$ is a $1$-cell
with endpoints in vertices $v_i$ and $v_{i+1}$.
Without loss of generality we may assume that
$v_i\in C$ implies $v_i=v$ or $v_i=w$.
Since $\bigcup\cA_{j-1}$ is closed and $v_m\not\in \cA_{j-1}$,
we see that also $e_m\not\in A_{j-1}$.
We will show that also $v_{m-1}\not \in A_{j-1}$.
If not, then
$e_m$ would have been placed on $L$ and then
removed together with  $v_{m}$ as a reduction pair
before the $j$th pass of the $\kwwhile$ loop.
Arguing by induction, we see that $v_0\not\in \cA_{j-1}$,
a contradiction.
\qed

\subsection{Quotients of reduction pairs.}
For a reduction pair $\alpha$
we write
$K_\alpha:=\cl \alpha_1$,
$\varphi_\alpha:=\varphi_{\alpha_1}$,
$\Dim{\alpha}:=\Dim{\alpha_1}$,
$K_\alpha^-:=\varphi_{\alpha}(\clb{B^n}_-)$,
$K_\alpha^+:=\varphi_{\alpha}(\clb{B^n}_+)$,
$K_\alpha^\setminus:=K\setminus\alpha_1\setminus\alpha_0$.
Set
\begin{eqnarray*}
    \cK_\alpha  &:=&\setof{\tau\in\cK\mid\tau\subset K_\alpha},\\
    \cK_\alpha^+&:=&\cK_\alpha\setminus\{\alpha_0,\alpha_1\}.
\end{eqnarray*}

\begin{lem}
\label{lem:cKalpha-subcomplexes}
We have $K_\alpha=\bigcup\cK_\alpha$ and $K_\alpha^+=\bigcup\cK_\alpha^+$.
In particular,  $K_\alpha$ and $K_\alpha^+$ are subcomplexes of $K$.
\end{lem}
\proof
Obviously $\bigcup\cK_\alpha\subset K_\alpha$.
To prove the opposite inclusion take an $x\in K_\alpha$.
If $x\in\alpha_1$, then $x\in\bigcup\cK_\alpha$,  because $\alpha_1\subset\cK_\alpha$.
Otherwise, $x\in \varphi_{\alpha_1}(S^{\Dim{\alpha}-1})$.
By assumption (i) of the definition of a regular facet we have $x\in\tau$ for some
cell $\tau\subset\varphi_{\alpha_1}(S^{\Dim{\alpha}-1})$. But $\varphi_{\alpha_1}(S^{\Dim{\alpha}-1})\subset\cK_\alpha$,
hence $x\in\bigcup\cK_\alpha$. This proves the first equality.

We have
\[
   K_\alpha=\varphi_{\alpha}(\clb{B^{\Dim{\alpha}}})
   =\varphi_{\alpha}(\opb{B^{\Dim{\alpha}}})\cup \varphi_{\alpha}(\clb{B^{\Dim{\alpha}}}_-)\cup \varphi_{\alpha}(\clb{B^{\Dim{\alpha}}}_+)
   =\alpha_1\cup\alpha_0\cup K_\alpha^+.
\]
But $\alpha_1\cap  K_\alpha^+=\emptyset$ by the definition of a cell and  $\alpha_0\cap  K_\alpha^+=\emptyset$
by the definition of the reduction pair. Therefore,
\[
  K_\alpha^+=K_\alpha\setminus\alpha_1\setminus\alpha_0=\bigcup\cK_\alpha^+,
\]
which shows that $K_\alpha$ is a subcomplex of $\cK$.
\qed

Let $\chi=(\chi_1,\chi_2): \clb{B^{\Dim{\alpha}}}\to\clb{B^{\Dim{\alpha}-1}}\times I$ be a homeomorphism such that
\begin{eqnarray*}
  \chi(\clb{B^{\Dim{\alpha}}}_-)&=&\clb{B^{\Dim{\alpha}-1}}\times\{1\}\cup S^{\Dim{\alpha}-1}\times I,\\
  \chi(\clb{B^{\Dim{\alpha}}}_+)&=&\clb{B^{\Dim{\alpha}-1}}\times\{0\}.
\end{eqnarray*}
Set $\chi_\alpha:=\varphi_\alpha\circ\chi^{-1}$.

Define the map $f_\alpha:K_\alpha\to K_\alpha^+$ by
\[
   f_{\alpha}(x):=\begin{cases}
      \chi_\alpha(\chi_1(\varphi_{\alpha}^{-1}(x)),0) & \text{ if $x\in\alpha_1\cup\alpha_0$,}\\
      x & \text{otherwise.}
   \end{cases}
\]
It is not difficult to check that $f_{\alpha}$ is well defined and continuous.

\begin{prop}
\label{prop:K-L-external-elementary-collapse}
The map $f_{\alpha}$ is a retraction of $K_\alpha$ onto $K_\alpha^+$, i.e. $f_{\alpha|K_\alpha^+}=\id_{K_\alpha^+}$.
Moreover, it is a homotopy inverse of the inclusion map $K_\alpha^+\subset K_\alpha$.
\end{prop}
\proof
The fact that $f_{\alpha}$ is a retraction is an elementary calculation.
To see that $f_{\alpha}$ is a homotopy equivalence,
define a homotopy $h:K_\alpha\times I\to K_\alpha^+$ by
\[
   h(x,t):=\begin{cases}
      \chi_\alpha(\chi_1(\varphi_{\alpha}^{-1}(x)),t\chi_2(\varphi_{\alpha}^{-1}(x))) & \text{ if $x\in\alpha_1\cup\alpha_0$,}\\
      x & \text{otherwise.}
   \end{cases}
\]
It is not difficult to check that $h$ is well defined, continuous and $h_0=f_{\alpha}$, $h_1=\id_{K_\alpha^+}$.
\qed

Define an equivalence relation $R_\alpha$ in $K$ by
\[
  R_\alpha:=\setof{(x,y)\in K^2\mid \text{$x=y$ or $x,y\in K_\alpha$ and $f_\alpha(x)=f_\alpha(y)$}}.
\]
Let $[x]_\alpha$ denote the equivalence class of $x$ with respect to $R_\alpha$ and
let $\kappa_\alpha:K\ni x\mapsto [x]_\alpha\in K/R_\alpha$ denote the quotient map.

\begin{prop}
\label{prop:KRalpha-Hausdorff}
  The space $K/R_\alpha$ is Hausdorff.
\end{prop}
\proof
Note that $[x]_\alpha=f_\alpha^{-1}(f_\alpha(x))$, therefore the equivalence classes of $R$ are compact.
Since $K$ is Hausdorff and $K_\alpha$ is closed in $K$,
in view of Proposition~\ref{prop:quotient-Hausdorff} and Proposition~\ref{prop:quotient-extension-usc}
it is enough to verify that the restriction of $R_\alpha$ to $K_\alpha$ is closed.
By Proposition~\ref{prop:quotient-usc} it suffices to show that for any $U$ open in $K_\alpha$
the set $U_R$ is open in $K_\alpha$. Assume this is not true.
Then, for some $x\in U_R$ no neighbourhood $V$ of $x$ satisfies $V\subset U_R$.
Thus,  Proposition~\ref{prop:cw-metrizable} allows us to select
a sequence $x_n\to x$ such that $[x_n]_R\not\subset U$. Let $y_n\in[x_n]\setminus U$.
Passing to a subsequence, if necessary, we may assume that $y_n\to y\in K\setminus U$.
Since $y_n\in[x_n]$, we have $f_\alpha(x_n)=f_\alpha(y_n)$. Thus, $f_\alpha(x)=f_\alpha(y)$.
It follows that $[x]_R=[y]_R$. In consequence $y\in U$, a contradiction.
\qed

Set $\kappa_{\alpha}^+:=\kappa_{\alpha\mid K_\alpha^+}$.

\begin{lem}
\label{lem:sub-alpha-quotient-bijection}
The restriction $\kappa_{\alpha}^\setminus:=\kappa_{\alpha\mid K_\alpha^\setminus}: K_\alpha^\setminus\to K/R_\alpha$
is a continuous bijection.
\end{lem}
\proof
The map $\kappa_{\alpha}^\setminus$ is continuous as a restriction of a continuous map.
To show that $\kappa_{\alpha}^\setminus$ is surjective,
take an $x\in K$. If $x\not\in K_\alpha$, then $\kappa_\alpha(x)=[x]_\alpha$.
Hence, assume $x\in K_\alpha$. Since $f_\alpha$ is a retraction onto $K_\alpha^+$,
we have $f_\alpha(x)=f_\alpha(f_\alpha(x))$. Thus $[x]_\alpha=[f_\alpha(x)]_\alpha=\kappa_{\alpha}^+(f_\alpha(x))$,
which shows that $\kappa_{\alpha}^+$ is surjective.
Assume in turn that
$[x_1]_\alpha=[x_2]_\alpha$  for some $x_1,x_2\in K$.
Then $f_\alpha(x_1)=f_\alpha(x_2)$.
If either $x_1$ or $x_2$ is not an element of $K_\alpha$, then both classes are singletons, hence $x_1=x_2$.
Thus, assume that $x_1,x_2\in K_\alpha^+$.
Since $f_\alpha$ is a retraction onto $K_\alpha^+$, we have
$x_1=f_\alpha(x_1)=f_\alpha(x_2)=x_2$, which shows that  $\kappa_{\alpha}^+$
is injective.
\qed

\begin{lem}
\label{lem:sub-alpha-quotient-homeo}
The restrictions
  $\kappa_{\alpha\mid K_\alpha^+}$ and
  $\kappa_{\alpha\mid K\setminus K_\alpha}$
are homeomorphisms onto their images.
\end{lem}
\proof
By Lemma~\ref{lem:sub-alpha-quotient-bijection} both maps are continuous bijections.
Since $K_\alpha$ is compact and $K/R_\alpha$ is Hausdorff, the first restriction is
a homeomorphism onto $K_\alpha^+$.
To see that the other restriction is a homeomorphism observe that
equivalence classes in $K\setminus K_\alpha$ are singletons, therefore
$\kappa_{\alpha}^{-1}(\kappa_{\alpha}(U))=U$ for any $U\subset K\setminus K_\alpha$.
In particular, for any
open set $U\subset K\setminus K_\alpha$ its image $\kappa_{\alpha}(U)$ is open,
which shows that also the other restriction is a homeomorphism.
\qed

\subsection{Collapses.}

We call the quotient space  $K/R_\alpha$
an {\em $\alpha$-collapse} of $K$ and write briefly $\kappa_\alpha$
for the quotient map $\kappa_{R_\alpha}$.

By performing a collapse we do not change the homotopy type and stay in the class
of CW complexes, as the following theorem shows.

\begin{thm}
\label{thm:alpha-quotient}
The quotient space $K/R_\alpha$ is a CW-complex with CW structure
$(\{\kappa_\alpha(\sigma)\}_{\sigma\in\cK\setminus\{\alpha_1,\alpha_0\}},
 (\{\kappa_\alpha\circ\varphi_\sigma\}_{\sigma\in\cK\setminus\{\alpha_1,\alpha_0\}})$.
Moreover, the complexes $K$ and $K/R_\alpha$ are homotopy equivalent.
\end{thm}
\proof
   First recall that by Proposition~\ref{prop:KRalpha-Hausdorff} $K/R_\alpha$ is a Hausdorff space.
   Let $\sigma\in\cK\setminus\{\alpha_1,\alpha_0\}$.
In order to prove property (CW0) of the definition
of CW complex it suffices to show that
$\kappa_\alpha(\sigma)$ is homeomorphic to $\sigma$.
By Lemma~\ref{lem:cKalpha-subcomplexes} the set $K_\alpha^+$ is a subcomplex, therefore
either $\sigma\subset K_\alpha^+$ or $\sigma\subset K\setminus K_\alpha$.
Thus, by Lemma~\ref{lem:sub-alpha-quotient-homeo},
in both cases $\kappa_\alpha(\sigma)$ is homeomorphic to $\sigma$,
as required. Properties (CW1) and (CW2) are straightforward.

By Proposition~\ref{prop:glueing} complex $K$ may be obtained from $K_\alpha$ by consecutively gluing a sequence
$(\sigma_j)$ of cells via the attaching maps $\theta_{\sigma_j}$.
It is not difficult to observe that applying to $K_\alpha/R_\alpha$ the same sequence of gluings
but with attaching maps  $\kappa_\alpha\circ\theta_{\sigma_j}$ we obtain the complex $K/R_\alpha$.
Therefore, the conclusion follows by induction from Theorem~\ref{thm:kozlov} and Theorem~\ref{thm:alpha-quotient}.
\qed

\begin{prop}
\label{prop:kappa-subcomplexes}
If $L$ is a subcomplex of $K$, then  $\kappa_\alpha(L)$ is a subcomplex of $K/R_\alpha$.
\qed
\end{prop}


Let $\cV$ be an acyclic vector field on $\cK$. Choose $\alpha_0,\alpha_1\in\cV$ such that
$\alpha_1=\cV(\alpha_0)$. Then $\alpha:=(\alpha_0,\alpha_1)$ is a reduction pair.
Define a partial map $\cV_{\alpha}:K/R_{\alpha}\to K/R_{\alpha}$ by
\begin{eqnarray*}
  \dom\cV_\alpha&:=&\setof{\kappa_\alpha(\sigma)\mid \sigma\in\dom\cV\setminus\{\alpha_0\}},\\
  \cV_\alpha(\kappa_\alpha(\sigma))&:=&\kappa_\alpha(\cV(\sigma)).
\end{eqnarray*}

\begin{thm}
\label{thm:recursive-acyclic-vector-field}
   If $\cV$ is acyclic, then  $\cV_{\alpha}$ is a well defined discrete vector field on $K/R_{\alpha}$
   and it is also acyclic.
\end{thm}
\proof
To show that $\cV_{\alpha}$ is a well defined discrete vector field on $K/R_{\alpha}$
consider $\beta_0\in\dom\cV\setminus\{\alpha_0\}$ and let $\beta_1:=\cV(\beta_0)$.
Then, $\beta:=(\beta_1,\beta_0)$ is a reduction pair in $K$. Moreover, $\beta\neq\alpha$.
Thus, $\beta_1,\beta_0\in K_\alpha^\setminus$. Since by Lemma~\ref{lem:sub-alpha-quotient-bijection}
$\kappa_\alpha$ is injective on $K_\alpha^\setminus$, we see that $\cV_{\alpha}$ is injective
and $\dom\cV_\alpha\,\cap\,\im\cV_\alpha=\emptyset$. Hence, we only need to show that
$\kappa_\alpha(\beta_0)$ is a regular facet of $\kappa_\alpha(\beta_1)$ in $\cK_\alpha$.
Note that (i) follows from Proposition~\ref{prop:kappa-subcomplexes}
and (ii) from Lemma~\ref{lem:sub-alpha-quotient-bijection}. To prove (iii) assume the contrary.
Then, $\kappa_\alpha(\varphi_\alpha(B^{\Dim{\alpha}}))\cap\kappa_\alpha(\beta_0)\neq\emptyset$.
Hence, we can choose an $x\in K_\beta^+$ and a $y\in\beta_0$
such that $x\neq y$ and $\kappa_\alpha(x)=\kappa_\alpha(y)$.
Since $\kappa_\alpha$ is bijective on $K_\alpha^\setminus$, we have $x,y\in K_\alpha$
and either $x$ or $y$ is not in $K_\alpha^+$. But $y\in\alpha_0\cup\alpha_1$ implies $\beta_0=\alpha_0$
or $\beta_0=\alpha_1$ and contradicts $\alpha\neq \beta$. Thus, $x\not\in K_\alpha^+$, i.e.
either $x\in\alpha_1$ or $x\in\alpha_0$. Consider $x\in\alpha_1$. Then $\alpha_1\cap K_\beta^+\neq\emptyset$
which implies $K_\alpha\subset K_\beta^+$ and
\[
  \emptyset\neq \beta_0\cap K_\alpha^+\subset \beta_0\cap K_\alpha^+\subset
  \beta_0\cap K_\alpha\subset  \beta_0\cap K_\beta^+,
\]
a contradiction. Consider $x\in\alpha_0$. Then, $\beta_0\cap K_\alpha^+$ and $\alpha_0\cap K_\beta^+$, i.e.
$\alpha_0$ is a facet of $\beta_1$ and $\beta_0$ is a facet of $\alpha_1$.
Thus, $\alpha_0,\alpha_1,\beta_0,\beta_1,\alpha_0$ is a cycle in $G_\cV$, a contradiction again.
This completes the proof that $\cV_{\alpha}$ is a well defined discrete vector field on $K/R_{\alpha}$.

We still need to prove that $\cV_{\alpha}$ is acyclic. Assume it is not.
Then, we have a cycle
\[
  \kappa_\alpha(\beta^1_0),\kappa_\alpha(\beta^1_1),\kappa_\alpha(\beta^2_0),\kappa_\alpha(\beta^2_1),\ldots
  \kappa_\alpha(\beta^n_0),\kappa_\alpha(\beta^n_1),  \kappa_\alpha(\beta^1_0)
\]
such that $\beta^i:=(\beta^i_0,\beta^i_1)\neq\alpha$,
$\kappa_\alpha(\beta^i_0))$ is a regular facet of $\kappa_\alpha(\beta^i_1)$
and $\kappa_\alpha(\beta^{i+1}_0))$ is a facet of $\kappa_\alpha(\beta^i_1)$.
Since $\beta_0^i,\beta_1^i\subset K_{\beta^i}$ and $K_{\beta^i}$ is a compact subset of $K_\alpha^\setminus$,
we see from Lemma~\ref{lem:sub-alpha-quotient-bijection} that $\kappa_{\alpha|K_{\beta^i}}$ is a homeomorphism,
hence $\beta^i_0$ is a regular facet of $\beta^i_1$.
If $\cl\beta_1^i\cap K_\alpha=\emptyset$ for all $i$, then
\begin{equation}
\label{eq:beta-sequence}
  \beta^1_0,\beta^1_1,\beta^2_0,\beta^2_1,\ldots  \beta^n_0,\beta^n_1,  \beta^1_0
\end{equation}
is a cycle in $G_\cV$, because $\kappa_\alpha$ restricted to $K\setminus K_\alpha$ is a homeomorphism onto
its image.
Hence, assume that $\beta^{i+1}_0$ is not a facet of $\beta^i_1$ for some $i$ but
$\kappa_\alpha(\beta^{i+1}_0)\subset \cl\kappa_\alpha(\beta^i_1)$. This is only possible if
$\alpha_0$ is a facet of $\beta^{i}_1$ and $\beta^{i+1}_0$ is a facet of $\alpha_1$.
Thus, we can insert $\alpha_0,\alpha_1$ between $\beta^{i+1}_0$ and $\beta^{i}_1$
in \eqref{eq:beta-sequence} obtaining a cycle in $G_\cV$, a contradiction.
\qed

A sequence $\mathbf{\alpha}=(\alpha^1,\alpha^2,\ldots,\alpha^m)$
of reduction pairs in $\cK$ is called {\em admissible}
if the associated collection $\{\alpha^1,\alpha^2,\ldots,\alpha^m\}$
of reduction pairs is a discrete vector field.
Obviously, given a discrete vector field $\cV$, any ordering
of its elements is admissible.
Assuming $\cV$ is acyclic and fixing an ordering $\mathbf{\alpha}$,
we can use Theorem~\ref{thm:recursive-acyclic-vector-field}
to recursively apply collapses of reduction pairs in
$\mathbf{\alpha}$. By Theorem~\ref{thm:alpha-quotient},
the resulting space is a CW complex
homotopy equivalent to the original complex $K$.
We call it the {\em Morse complex} of the vector field $\cV$
and denote $K/\cV$. This is some abuse of language, because
the construction depends on the ordering $\mathbf{\alpha}$
but we assume the ordering is implicit.
Summarizing, we get the following fundamental result of discrete Morse theory.

\begin{cor}
\label{cor:fundamental-DMT}
Let $K$ be a finite CW complex and $\cV$ an acyclic discrete vector field
on the CW structure $\cK$ of $K$. Then, the complexes $K$ and $K/\cV$
are homotopy equivalent.
\end{cor}
\qed

\section{Reductions.}
\label{sec:reductions}
\subsection{External collapses}
We say that the reduction pair $\alpha$ is {\em free} if
$K_\alpha^\setminus$ is a subcomplex of $K$.

\begin{thm}
\label{thm:external-collapse}
If $\alpha$ is free then the complexes $K$ and $K_\alpha^\setminus$ are homotopy equivalent.
\end{thm}
\proof
We know from Lemma~\ref{lem:sub-alpha-quotient-bijection} that the restriction $\kappa_{\alpha|K_\alpha^\setminus}$
is a continuous bijection onto $K/R_\alpha$.
Thus, it is a homeomorphism, because $K_\alpha^\setminus$, as a subcomplex, is compact.
Hence, the conclusion follows from Theorem~\ref{thm:alpha-quotient}.
\qed

When $\alpha$ is free, we say that  $K_\alpha^\setminus$ is
an {\em external collapse} of $K$. Otherwise we say that $K/R_\alpha$ is an {\em internal collapse} of $K$.
On the algorithmic side external collapses are advantageous,
because to compute them it is enough to remove the reduction pair
from the list $\cK$ of cells of $K$. Thus, it is convenient to make as many external collapses as possible.
We say that $K$ collapses onto the subcomplex $L$ if there exists a sequence of
subcomplexes $L=L_0\subset L_1\subset \ldots \subset L_m=K$ such that
$L_{j-1}$ is an external collapse of $L_j$.
Note that for convenience we admit
the case when $m=0$, i.e. $L=K$.

An immediate consequence of Proposition~\ref{prop:K-L-external-elementary-collapse} is the following
proposition.
\begin{prop}
\label{prop:K-L-external-collapse}
If $K$ collapses onto a subcomplex $L$ then the inclusion $L\subset K$ is a homotopy equivalence. \qed
\end{prop}

\subsection{Collapsible subcomplexes.}
Complex $K$ is called {\em collapsible} if it collapses onto a vertex.

\begin{thm}
\label{thm:collapsible-subset}
Assume $L$ is a collapsible subcomplex of $K$.
Then, $K$ and $K/L$ are homotopy equivalent.
\end{thm}
\proof
Let $L_m\subset L_{m-1}\subset \ldots \subset L_0=L$ be a sequence of subcomplexes of $L$ such that
$L_{j}$ is an external collapse of $L_{j-1}$ by a reduction pair $\alpha^j$ and $L_m=\{*\}$ consists only of a vertex.
Let $\kappa_j:L_j\to L_j/R_{\alpha^j}$ be the quotient map and let
$\kappa_*:=\kappa^m\circ\kappa^{m-1}\circ\cdots\kappa^1$.
By Proposition~\ref{prop:sub-reduction-pair}
each $\alpha^j$ is a reduction pair in $K$.
Thus, by Theorem~\ref{thm:external-collapse} $K$ is homotopy equivalent to
\[
  K/R_{\alpha^1}/R_{\alpha^2} \cdots /R_{\alpha^m}
\]
and by Proposition~\ref{prop:quotient-iterated} this space is homeomorphic to $K/R_{*}$,
where $R_{*}$ is a relation on $X$ defined by $xR_* y$ if and only if $\kappa_*(x)=\kappa_*(y)$.
But  $\kappa_*$ is constant on $L$ and a homeomorphism on $K\setminus L$.
Hence, $K/R_{*}$ and $K/L$ are homeomorphic, and consequently $K$ and $K/L$ are homotopy equivalent.
\qed

\begin{thm}
\label{thm:collapsability-onto-via-dvf}
Assume $K$ is a regular CW complex
and $L$ is a subcomplex with CW structure
$\cL\subset \cK$. Then, $K$ collapses onto $L$ if and only if there
exists an acyclic discrete vector field on $\cK$ such that $\crit\cV=\cL$.
\end{thm}
\proof
   Assume first $K$ collapses onto $L$. Then, we have
   subcomplexes $L=L_0\subset L_1\subset \ldots \subset L_m=K$ such that
   $L_{j-1}$ is an external collapse of $L_j$ resulting from
   removing a reduction pair $\alpha^j$. It is straightforward
   to verify that the collection$\{\alpha^1,\alpha^2,\ldots\alpha^n\}$ is an acyclic discrete vector field on $\cK$ whose set of critical cells is precisely $\cL$.

   To prove the opposite implication assume $\cV$ is an acyclic
   discrete vector field on $\cK$ such that $\crit \cV=\cL$.
   Let $n$ be the cardinality of $\cK\setminus\cL$. We proceed by
   induction on $n$. For $n=0$ the conclusion is trivial.
   Thus, consider the case $n>0$. Then, we can choose a non-critical cell in $\sigma\in\cK$.
   Since $\cV$ is acyclic, we can choose a cell $\tau$ which is inferior to $\sigma$,
   that is $\tau$ is a minimal cell less than or equal to $\sigma$ in the partial order induced by $\cV$.
   We claim that $\tau$ is not a critical cell of $\cV$. To see this
   consider a path $\tau=\tau_0,\tau_1,\ldots,\tau_m=\sigma$ in $G_\cV$.
   Since $\tau$ is  inferior to $\sigma$, such a path obviously exists.
   If $\tau$ is critical, then we consider the smallest index $k$
   such that $\tau_{k-1}$ is critical and $\tau_{k}$ is not critical.
   Since the partial order induced by $G_\cV$ goes against the facet partial order
   only for the vectors of $\cV$, we have $\tau_k\in\dom\cV$
   and $\tau_k\subset\cl\tau\subset L$, because $\tau$, as a critical cell is a subset
   of $L$ and $L$ as a subcomplex is closed. It follows that $\tau_k$ is critical, a contradiction. Therefore $\tau\in\dom\cV$. Let $K':=K\setminus\tau\setminus\cV(\tau)$.
   To prove that $K'$ is closed, assume the contrary.
   Then, there exists an $x\in\cl K'\setminus K'$. Hence, $x\in\tau\cap\cV(\tau)$
   and $x\in\cl\sigma'$ for some cell $\sigma'\subset K'$.
   By Proposition~\ref{prop:intersection-implies-face}  $\tau$ is a face of $\sigma'$
   and by Proposition~\ref{prop:facets} $\tau$ is a facet of a face $\sigma''$ of $\sigma'$. Thus $\sigma''$ precedes $\tau$ in the partial order induced by $G_\cV$,
   and consequently $\tau$ is not  inferior to $\sigma$. This contradiction shows
   that $K'$ is a subcomplex of $K$ and consequently $K$ collapses onto $K'$.
   By induction assumption $K'$ collapses onto $L$. Therefore, $K$ collapses onto $L$.
\qed

As an immediate consequence of Theorem~\ref{thm:collapsability-onto-via-dvf}
we obtain the following corollary.
\begin{cor}
\label{cor:collapsability-via-dvf}
A regular CW complex is collapsible if and only if it admits a discrete vector field
whose set of critical cells consists of precisely one vertex.
\qed
\end{cor}

\subsection{Redundant configurations.}
Let $K$ be a regular CW complex, $L$ its subcomplex and  $\cK$ and $\cL$
the corresponding CW structures. We refer to the pair $(K,L)$ as a {\em configuration}.
We say that the configuration $(K,L)$ is {\em redundant} if there exists a discrete vector field
$\cV$ on $\cK$ such that $\crit \cV=\cL$.
We say that subcomplexes $L$, $M$ of $K$ form a {\em redundant decomposition}
of $K$ if $K=L\cup M$ and the configuration $(M,L\cap M)$ is redundant.

\begin{thm}
\label{thm:redundancy}
Assume subcomplexes $L$ and $M$ of a regular CW complex $K$ form
a redundant decomposition. Then, the inclusion $L\subset K$ is a homotopy equivalence.
\end{thm}
\proof
Let $\cV$ be a discrete vector field on $M$ such that $\crit \cV$ consists of the cells
in $L\cap M$. Extend it to all cells in $K$ by making all cells not in the domain of $\cV$ critical. It follows from Theorem~\ref{thm:collapsability-onto-via-dvf}
that $K$ collapses onto $L$, thus we get the assertion from
Proposition~\ref{prop:K-L-external-collapse}.
\qed

The advantage of Theorem~\ref{thm:redundancy} is twofold. We may use it to reduce
the complex without changing its homotopy type by removing a subcomplex from
a redundant decomposition. We may also simplify the complex
by quotienting out a maximal collapsible subcomplex. We can  construct
such a subcomplex $L\subset K$ by starting from a vertex and recursively
adding to $L$ a subcomplex $M$ such
that $L$ and $M$ form a redundant decomposition of $L\cup M$.
However, to make it work in practice we need a quick method to test for
redundant configurations.

This is possible in the setting of tessellations discussed next.

\subsection{Tesselated spaces and lattice complexes}

We say that a CW structure $\cK$ on a pure, regular CW complex $X$, finite or not,
is a {\em tessellation} if the closures of any two top-dimensional cells $\sigma_1$, $\sigma_2$
in $\cK$ are isomorphic, i.e. there exists
a homeomorphism $h:\cl \sigma_1\to\cl\sigma_2$ mapping cells onto
cells.
A CW complex is a {\em tesselated CW complex} if its CW structure is a tessellation.
Obviously, every pure subcomplex of a tesselated space is a tesselated space.

A special way to define a tessellation is to consider a {\em lattice} $L \subseteq \RR^n$,
i.e. an additive subgroup of $\RR^n$
generated by some choice of $n$ linearly independent vectors
$V_L= \{v_1, \cdots, v_n\}$. Any $v\in L$ determines a Dirichlet--Voronoi cell
\[
   D_L(v) := \{x \in \RR^n \mid~~ ||v-x|| \leq ||w-x|| ~ \forall w\in L\}.
\]
For $x\in\RR^n$ set
\[
   \chi_L(x):=\setof{v\in L\mid x\in D_L(v)}
\]
and define an equivalence relation $R\subset\RR^n\times\RR^n$ by
$xRy$ if and only if $\chi_L(x)=\chi_L(y)$. One can show that the equivalence classes
of this relation provide a CW structure of $\RR^n$ which is a tessellation.
We  call it a {\em lattice tessellation} and denote it $\cK_L$.
When $V_L$ consists of $n$ orthogonal vectors,
the sets $D_L(v)$ are $n$-dimensional hypercubes
and the tessellation is the classical cubical decomposition of $\RR^n$.
By identifying $\RR^n$ with the hyperplane
\[
\{(x_1, \ldots, x_{n+1}) \in \RR^{n+1} \ : \ x_1 + \cdots +x_{n+1}=0\}
\]
and taking as  $V_L$ the set
$v_1=(1-n,1,1, \ldots,1 , 1)$, $v_2=(1,1-n,1, \ldots,1, 1)$, \ldots,
$v_n=(1,1,1,\ldots, 1-n,1)$ we obtain a tessellation of $\RR^n$ consisting of
$n$-dimensional permutahedra.

A pure, finite CW complex whose CW structure is a subfamily of a lattice tessellation
will be called a {\em lattice CW complex}.
Lattice CW complexes may be easily and efficiently stored in bitmaps.

Fix a lattice $L$.
Two top-dimensional cells in $\cK_L$ are called {\em neighbours} if
the intersection of their closures is non-empty.
Consider a top-dimensional cell $\sigma\in\cK_L$
and a lattice CW complex $K$ with its CW structure $\cK\subset\cK_L$.

Observe that $K\setminus\sigma$ and $\cl\sigma$ are
also lattice CW complexes. Their intersection is a subcomplex of $X$.
We will call it the {\em contact complex} of $\sigma$ and $K$. Note that this makes
sense both if $\sigma\subset K$ and $\sigma\cap K=\emptyset$.

Let $N_K(\sigma)$ denote the set of neighbours of $\sigma$ contained in $K$.
Then
\[
  (K\setminus\sigma)\cap\cl\sigma=\bigcup\setof{\cl\sigma\cap\cl\tau\mid\tau\in N_K(\sigma), \tau\neq\sigma},
\]
i.e. the contact space is uniquely determined by the neighbours of $\sigma$ in $K$.

The toplex $\cl\sigma$ is said to be {\em redundant} in $K$ if the configuration
$(\cl\sigma,K\setminus\sigma\cap\cl\sigma)$ is redundant.
If $\sigma$ has $k$ neighbours, then there are $2^k$ configurations to test for redundancy.
Since all toplexes in a lattice tessellation are isomorphic, tests done for one
toplex may be reused for other toplexes.
Even better, for reasonably small number of neighbours
we can test all configurations of a toplex
for redundancy and store them in a lookup table for later quick redundancy tests.
The number of neighbours for a cubical tessellation in $\RR^n$ is $3^n-1$.
The number of neighbours for a permutahedral tessellation in $\RR^n$ is $2^{n+1}-2$.
Thus, the lookup tables  may be computed and stored in less than 1GB
of memory for cubical tessellations up to dimension
$3$, requiring $2^{26}$ bits or $2^{23}$ bytes in dimension $3$ and for permutahedral tessellations
up to dimension $4$, requiring $2^{30}$ bits or $2^{27}$ bytes in dimension $4$.

\section{Reduction algorithms.}
\label{sec:red-algo}

\subsection{Shaving and collapsing.}
Two reduction algorithms based on redundancy tests are presented
in Tables~\ref{tab:shaving} and~\ref{tab:collapsibleSubset}.
The first, often called shaving, removes redundant toplexes as long as they exist.
The second constructs a possibly large collapsible subcomplex.

\begin{table}
\begin{tabular}{|c|}
\hline
~~\\
\begin{minipage}{120mm}
\begin{algo}
\label{alg:shaving}
\upshape shaving(\texttt{collection of toplexes} $K$)\\
\0{$Q:=\emptyset$;} \\
\0{\kwwhile there exists a redundant toplex $\sigma$ in $K$ \kwdo} \\
\1{  remove $\sigma$ from $K$;}\\
\1{  enqueue $N_K(\sigma)$ in $Q$;}\\
\1{  \kwwhile $Q\neq\emptyset$ \kwdo} \\
\2{    $\sigma$=dequeue$(Q)$;}\\
\2{    \kwif $\sigma$ is a redundant in $K$ \kwdo} \\
\3{      remove $\sigma$ from $K$;}\\
\3{      enqueue $N_K(\sigma)$ in $Q$;}\\
\0{\kwreturn K;}
\end{algo}
\end{minipage}\\
~~\\
\hline
\end{tabular}
\vspace*{3mm}
\caption{Shaving}
\label{tab:shaving}
\end{table}

\begin{thm}
\label{thm:shaving-algo}
Algorithm~\ref{alg:shaving} applied to a collection of toplexes representing
a lattice space always stops and returns a homotopy equivalent
subcomplex with no redundant toplexes.
\end{thm}
\proof
The conclusion follows immediately from Theorem~\ref{thm:redundancy}.
\qed

\begin{table}
\begin{tabular}{|c|}
\hline
~~\\
\begin{minipage}{120mm}
\begin{algo}\upshape
\label{alg:collapsibleSubset}
collapsibleSubset(\texttt{collection of toplexes} $K$)\\
\0{$C:=Q:=\emptyset$;} \\
\0{choose any toplex $\sigma$ in $K$ and insert it into $C$;} \\
\0{enqueue $N_{K\setminus C}(\sigma)$ in $Q$;}\\
\0{\kwwhile $Q\neq\emptyset$ \kwdo} \\
\1{  $\sigma$:=dequeue$(Q)$;}\\
\1{  \kwif $\sigma$ is redundant in $C$ \kwdo} \\
\2{    insert $\sigma$ into  $C$;}\\
\2{    enqueue $N_{K\setminus C}(\sigma)$ in $Q$;}\\
\0{\kwreturn C;}
\end{algo}
\end{minipage}\\
~~\\
\hline
\end{tabular}
\vspace*{3mm}
\caption{Collapsible subset construction}
\label{tab:collapsibleSubset}
\end{table}

\begin{thm}
\label{thm:collapsible-subset-algo}
Algorithm~\ref{alg:collapsibleSubset} applied to a collection of toplexes representing
a lattice space always stops and returns  a collapsible subcomplex.
The resulting quotient complex is homotopy equivalent to the complex on input.
\end{thm}
\proof
The conclusion follows immediately from Theorem~\ref{thm:redundancy} and
Theorem~\ref{thm:collapsible-subset}.
\qed

\subsection{$C$-structures.}

In order to obtain an algorithm based on acyclic vector fields and
Theorem~\ref{thm:edge-group} it is convenient to
consider the following concept.

Let $C$ be a finite set decomposing into the disjoint union $C=C_0\cup C_1\cup C_2$.
By a {\em $C$-structure} we mean a quadruple $(C,e_{-1},e_{+1},d)$ such that
$(C_0,C_1,e_{-1},e_{+1})$ is a directed multigraph and $d:C_2\to\Cycles(C_1)$, where
$\Cycles(C_1)$ stands for the set of cycles in $C_1$.

We write $\Dim{c}=j$ if and only if $c\in C_j$.

For $c\in C$ we define its {\em boundary} by
\[
  \bd c:=\begin{cases}
           \emptyset & \text{ if $\Dim{c}=0$,}\\
           \{e_{-1}(c),e_{+1}(c)\} & \text{ if $\Dim{c}=1$,}\\
           \setof{e_1,e_2,\ldots e_n} & \text{ if $\Dim{c}=2$ and $d(c)=e^{\epsilon_1}_1e^{\epsilon_2}_2\cdots e^{\epsilon_n}_n$.}
         \end{cases}
\]

With every $C$-structure $C$ we associate a $2$-dimensional CW complex $K$
defined as follows. For each $c\in C$
we have a formal cell $\cell{c}:=\opb{B^{\Dim{c}}}\times\{c\}$.
However, to simplify notation, in the sequel we identify $\cell{c}$ with $c$.
Let
\[
K^j:=\bigcup\setof{c\in C\mid \;\Dim{c}\leq j}.
\]
For $c\in C_1$ and $t\in \clb{B^1}$ we define the characteristic map
$\varphi_c:\clb{B^1}\to K^1 $ by
\[
  \varphi_c(t):=\begin{cases}
                    (t,c) & \text{if $t\in \opb{B^1}$,}\\
                    e_{t}(c) & \text{otherwise.}
                 \end{cases}
\]

To define the characteristic map  $\varphi_c:\clb{B^2}\to K^2 $
for a $c\in C_2$ assume $d(c)=e^{\epsilon_1}_1e^{\epsilon_2}_2\cdots e^{\epsilon_n}_n$
and take an $x\in\clb{B^2}$. Set
\[
  \varphi_c(x):=\begin{cases}
                    (x,c) & \text{if $x\in \opb{B^2}$,}\\
                    \varphi_{e_k}^{\epsilon_k}(\frac{nt}{2\pi}-k) & \text{if $x=e^{2\pi ti}$ and
                    $\frac{2\pi}{n} k\leq t\leq \frac{2\pi}{n}(k+1)$.}
                 \end{cases}
\]

It is not difficult to verify that the cells $c$ together with the characteristic
maps $\varphi_c$ define a CW complex. We denote it by $\CW(C)$.
Conversely, given a finite, regular, $2$-dimensional CW complex $K$ it is straightforward
to obtain a $C$-structure $C$ such that $\CW(C)$ is isomorphic to $K$.

Given a reduction pair $(\alpha_0,\alpha_1)$
in $\CW(C)$ we define a $C$-structure $C'$ as follows.
Let $C':=C\setminus\{\alpha_0,\alpha_1\}$.
To define $e'_\epsilon:C'_1\to C'_0$ and $d':C'_2\to \Cycles(C'_1)$,
consider first the case $\Dim{\alpha_1}=1$. By the definition of a reduction pair
$\alpha_0$ is a regular facet of $\alpha_1$.
This implies that $\bd \alpha_1$ consists of precisely two
elements one of which is $\alpha_0$. Let $\alpha_0^*$ be the other.
Set
\begin{eqnarray*}
e'_\epsilon(c)&:=&\begin{cases}
                 \alpha_0^* & \text{if $e_{\epsilon}(c)=\alpha_0$,}\\
                 e_\epsilon(c) & \text{otherwise,}\\
               \end{cases}\\
         d'(c)&:=&d(c).
\end{eqnarray*}
Consider the case $\Dim{\alpha_1}=2$.
Since $\alpha_0$ is a regular facet of $\alpha_1$,
the facet $\alpha_0$ appears precisely once in the cycle
$d(\alpha_1)=e^{\epsilon_1}_1e^{\epsilon_2}_2\cdots e^{\epsilon_n}_n$.
Without loss of generality we may assume that $\alpha_0=e_1$. Let
$\alpha_0^*:=\left(e^{\epsilon_2}_2e^{\epsilon_3}_3\cdots e^{\epsilon_n}_n\right)^{-\epsilon_1}$.
Then $\alpha_0^*$ and $\alpha_0$ have the same initial and terminal vertex.
Moreover, $\alpha_0$ does not appear in $\alpha_0^*$.
Set
\[
d'(c):=S_{\alpha_0,\alpha_0^*}(d(c)),
\]
where $S_{\alpha_0,\alpha_0^*}$ denotes the $\alpha_0,\alpha_0^*$
substitution in $d(c)$.

It is not difficult to verify  that $(C',e'_{-1},e'_{+1},d')$ is a $C$-structure. We call it the {\em $\alpha$-collapse} of $C$.
One easily verifies the following proposition.
\begin{prop}
\label{prop:alpha-collapse}
The CW complex $\CW(C')$ associated with the $\alpha$-collapse $C'$
is homeomorphic to the quotient complex $\CW(C)/R_\alpha$.
\qed
\end{prop}

Algorithm~\ref{alg:fundamentalGroup} computing a presentation
of the homotopy group of a regular CW complex $K$
represented as a collection of top-dimensional cells is presented
in Table~\ref{tab:fundamentalGroup}.

\begin{table}
\begin{tabular}{|c|}
\hline~~\\
\begin{minipage}{120mm}
\begin{algo}\upshape
\label{alg:fundamentalGroup}
fundGroup(\texttt{collection of toplexes} $K$)\\
\0{$K:=$ shaving$(K)$;} \\
\0{$A:=$ collapsibleSubset$(K)$}; \\
\0{$C:=$ C-structure of $K/A$;}\\
\0{$V:=$ discreteVectorField$(C)$;} \\
\0{\kwfor\kweach $\alpha\in V$ \kwdo}\\
\1{ assign to $C$ the $\alpha$-collapse of $C$;}\\
\0{\kwreturn $(C_1,d(C_2))$;}\\
\end{algo}
\end{minipage}\\
~~\\
\hline
\end{tabular}
\vspace*{3mm}
\caption{Fundamental group presentation of a regular CW complex.}
\label{tab:fundamentalGroup}
\end{table}

\begin{thm}
\label{thm:fundamentalGroup-algo}
Assume Algorithm~\ref{alg:fundamentalGroup}
is called with a collection of top cells of a CW structure of
a lattice space $K$. It always stops and returns a presentation,
up to an isomorphism, of the fundamental group of $K$.
\end{thm}
\proof
The theorem follows immediately from Theorem~\ref{thm:shaving-algo},
Theorem~\ref{thm:collapsible-subset-algo},
Theorem~\ref{thm:discreteVectorField-algo},
 Theorem~\ref{thm:alpha-quotient}, Proposition~\ref{prop:alpha-collapse},
 Theorem~\ref{thm:recursive-acyclic-vector-field}
and Theorem~\ref{thm:edge-group}.
\qed

\section{Application to knot classification.}
\label{sec:applications}

\subsection{Knot classification.}
Recall  that two knots $K,K'\subset\RR^3$ are {\em equivalent } \cite{Moran-1983,GordonLuecke-1989} if there exists
a homeomorphism  $f:\RR^3\to\RR^3$ such that $K'=f(K)$.
They are {\em orientation equivalent} \cite{Moran-1983} or {\em isotopic} \cite{GordonLuecke-1989}
if the homeomorphism is an {\em ambient isotopy} of knots, that is it is isotopic to the identity.
Gordon and Luecke \cite[Theorem 1]{GordonLuecke-1989} prove that
$K,K'$ are equivalent if and only if the complements $\RR^3\setminus K$ and
$\RR^3\setminus K'$ are homeomorphic.
Thus, topological invariants of the complements of knots
may be used to distinguish between non-equivalent knots.
They cannot be used to distinguish between non-isotopic knots.
But, by a result of Fisher \cite{Fi60} two equivalent but non-isotopic knots
differ only by the mirror image.

Let $T$ be an invariant of topological spaces,
that is $T(X)=T(Y)$ if $X$ and $Y$ are homeomorphic.
Then, obviously $T(\RR^3\setminus K)=T(\RR^3\setminus K')$ if $K$ and $K'$
are equivalent.
We say that $T$ {\em classifies a family of knots $\cK$} if for any two non-equivalent knots $K,K'\in\cK$ we have $T(\RR^3\setminus K)\neq T(\RR^3\setminus K')$.
Note that if the family
$\cK$ is finite and the invariant $T$ is algorithmically computable, then,
at least theoretically, one can verify if $T$ classifies $\cK$ in finite computation.
A natural candidate for a good $T$ is the fundamental group.
Unfortunately, by the results of Novikov \cite{Novikov-1955} and Boone \cite{Boone-1959},
in general it is not possible to decide algorithmically from finite presentations of two groups if they are non-isomorphic. However, this finite presentation may be used to compute
some further, algorithmically decidable, invariants. This is the path we used
in \cite{Brendel-etal-2014}. Recall that
given a knot $K \subset \RR^3$, using Algorithm~\ref{alg:fundamentalGroup}
we can compute a presentation of the fundamental group
$G_K:=\pi(D^3\setminus M_K)$, where $D^3$
is a cubical ball in $\RR^3$ containing the knot $K$ and $M_K$
is a cubical neighbourhood of $K$ containing $K$ as a deformation retract.
Using this presentation and the GAP software library, we can compute some
algebraic invariants of $G_K$. In particular,
given integers $n,c\ge 1$, $m\ge 0$ we study in \cite{Brendel-etal-2014} the invariant $I^{[n,c,m]}(K)$.
It is defined by
\[
I^{[n,c,m]}(K) := \{H_m(S/\gamma_{c+1}S,\mathbb Z)\ : \ S\le G_K, |G_K:S|\le n\},
\]
where $H_m(G,\mathbb Z)$ denotes the integral homology of group $G$ and $\gamma_k G$
is defined recursively by $\gamma_1 G:=G$, $\gamma_{k+1} G:=[\gamma_k G,G]$.
Trying to balance the computational expenses between the cost of computing the invariant
$I^{[n,c,m]}(K)$ and its ability to distinguish between knots we chose the invariant
$I^n(K):=I^{[n,1,1]}(K)$ with $c=m=1$ and varying $n$.
Using a GAP implementation of Algorithm~\ref{alg:fundamentalGroup}
we  provide in \cite{Brendel-etal-2014}
a computer assisted proof of the following result.

\begin{thm}
\label{thm:cross-11}
The invariant $I^6(K)$ distinguishes, up to mirror image,
between ambient isotopy classes of all prime knots that admit
planar diagrams with eleven or fewer crossings.
\qed
\end{thm}

\subsection{Classifying index.}
Using a C++  implementation of Algorithm~\ref{alg:fundamentalGroup}
and an improved strategy of computing the classifying invariant  we are able
to strengthen Theorem~\ref{thm:cross-11} (see Theorem~\ref{thm:cross-14} below).
The computational challenge of this extension lies, in particular, in the exponential growth
of the number $N(c)$ of prime knots with  $c$ crossings. Prime knots up to $16$ crossing
were tabulated by  Hoste, Thistlethwaite and Weeks in \cite{Hoste-etal-1998}.
The number of such knots for a given $n$ is presented in Table~\ref{tab:prime-knots-growth}.
In our computations we used the Hoste-Thistlethwaite tabulation available from \cite{KnotAtlas}.

\begin{table}
\begin{tabular}{||c ||c |c |c |c |c |c |c ||}
  \hline
  \hline
  $c$ & 3 & 4 & 5 & 6 & 7 & 8 & 9  \\
  \hline
  $N(c)$ & 1 & 1 & 2 & 3 & 7 & 21 & 49 \\
  \hline
  \hline
  $c$ & 10 & 11 & 12 & 13 & 14 & 15 & 16 \\
  \hline
  $N(c)$ & 165 & 552 & 2176 & 9988 & 46972 & 253293 & 1388705 \\
  \hline
  \hline
\end{tabular}
\vspace*{3mm}
\caption{The growth of the number of prime knots $N(c)$ with the number of crossings $c$.}
\label{tab:prime-knots-growth}
\end{table}

There are only $801$ prime knots up to 11 crossings but $59937$ prime knots
up to 14 crossings. The additional challenge is in the cost of computing the invariant
$I^n$, particularly when $n$ grows. An observation we made is that for most knots relatively low values
of $n$ suffice to distinguish them from other knots. But, some knots require relatively high values of $n$.
Thus, in order to decrease the computational cost, it makes sense to take as the classifying invariant the pair $(n,I^n(K))$
for a possibly low $n$ selected individually for each knot.
Assuming the invariant classifies a family of knots $\cK$, given
a knot $K\in\cK$ we refer to $n$ in the invariant $(n,I^n(K))$ as the {\em classifying index
of $K$}. By the {\em classifying index of the family $\cK$} we mean the maximum of all classifying indexes
of individual knots  $K\in\cK$.
The idea of classifying indexes  leads to Algorithm~\ref{alg:indexing}
presented in Table~\ref{tab:indexing}.
Note that the $\kwwhile$ loop in Algorithm~\ref{alg:indexing} is potentially
infinite. But, we have the following theorem.

\begin{table}
\begin{tabular}{|c|}
\hline
~~\\
\begin{minipage}{120mm}
\begin{algo}\upshape
\label{alg:indexing}
classify(\texttt{collection of knots} $\cK$)\\
\0{$Q:=$ empty queue; // to process $(K,n)$ }\\
\0{$A:=$ empty set; // $(K,n)$ added to $Q$}\\
\0{$I:=$ empty dictionary; // $I[i]$ is a list of $(K, n)$, s.t. $i=I^n(K)$ }\\
\0{$C:=$ empty dictionary; // $C[K]$ final invariant for $K$ } \\
\0{$D:=$ empty dictionary; // $D[i]=K$ iff $C[K]=i$} \\
\0{$N:=0$; // classifying index}\\
\0{\kwforeach $K\in \cK$ \kwdo}\\
\1{enqueue  $(K, 2)$ in $Q$;}\\
\1{insert  $(K, 2)$ into $A$;}\\
\0{\kwwhile $Q$ is non-empty \kwdo} \\
\1{  $(K, n)$:=dequeue$(Q)$;}\\
\1{  $N:=\texttt{max}(N,n)$;}\\
\1{  $i:=I^n(K)$;}\\
\1{  \kwif $i$ is not a key in $I$ \kwthen}\\
\2{     $I[i] :=$ empt list;}\\
\2{     append $(K,n)$ to list $I[i]$;}\\
\2{     $D[i] := K$;}\\
\2{     $C[K]:=(n,i)$;}\\
\1{  \kwelse}\\
\2{     append $(K,n)$ to list $I[i]$;}\\
\2{     remove key $i$ from $D$;}\\
\2{     \kwforeach $(\bar{K},\bar{n}) \in I[i]$ \kwdo}\\
(*)\2{        remove $\bar{K}$ from $C$;}\\
\3{        \kwif $(\bar{K},\bar{n} + 1) \notin A$ \kwthen}\\
(**)\3{            \!\!enqueue $(\bar{K},\bar{n} + 1)$ in $Q$;}\\
\4{            insert $(\bar{K},\bar{n} + 1)$ into $A$;}\\
\0{\kwreturn (C,D,N);}
\end{algo}
\end{minipage}\\
~~\\
\hline
\end{tabular}
\vspace*{3mm}
\caption{Classification algorithm.}
\label{tab:indexing}
\end{table}

\begin{thm}
\label{thm:indexing}
Assume Algorithm~\ref{alg:indexing} starts with a collection of knots
$\cK$ on input and stops.
Then, the triple $(C,D,N)$ on output, consisting of dictionaries $C$, $D$
and an integer $N$, has the following features.
\begin{itemize}
   \item[(a)]   The integer $N$ is the classifying index of $\cK$.
   \item[(b)]  The keys of $C$ are all knots in $\cK$.
   \item[(c)] The keys of $D$  constitute a collection of invariants $I^n(K)$,
                one per each knot $K\in\cK$. Given an invariant $i$ in the keys of $D$,
                the associated value $D[i]$ is the unique knot $K$ in $\cK$ such that $C[K]=(n,I^n(K))=(n, i)$.
                Then, the number $n$ is the classifying index of $K$.
   \item[(d)]  Both dictionaries $C$ and $D$ are injective.
                In particular, for any two knots $K,L\in\cK$ the invariants $C[K]$ and $C[L]$ are different.
\end{itemize}
\end{thm}
\proof
Given a variable $x$ in Algorithm~\ref{alg:indexing} we let $x^{(l)}$
denote the state of  $x$  after $l$ iterations of the $\kwwhile$ loop.
In particular, for $l=0$ this is the state just before entering the $\kwwhile$ loop.
We will show by induction that for every $l\geq 0$ we have the following properties:
\begin{enumerate}[label=(\roman*)]
\item \label{item:alg-a} $A^{(l)} = \bigcup_{k=0}^{l} Q^{(k)}$,
\item \label{item:alg-seq} the sequence $\setof{n^{(i)}}_{i=1}^l$ is not descending,
\item \label{item:alg-uniq} $(K^{(l)}, n^{(l)}) \ne (K^{(l')}, n^{(l')})$ for each $l' < l$,
\item \label{item:alg-C} for each $K\in\cK$ the value of $C^{(l)}[K]$ is set if and only if $K\in\cK\setminus\cL^{(l)}$, where
\[
  \cL^{(l)} := \setof{K\in\cK \mid\exists_{n\in{\scriptsize\NN}} (K,n) \in Q^{(l)}},
\]
\item \label{item:alg-inj} $C^{(l)}$ and $D^{(l)}$ are injective,
\item \label{item:alg-value} the value $D^{(l)}[i]$ is the unique $K\in\cK$ such that
\begin{equation}
\label{eq:alg-value}
  C^{(l)}[K]=(n,I^n(K))=i.
\end{equation}
\end{enumerate}
First observe that  dictionaries $C^{(0)}$, $D^{(0)}$, $I^{(0)}$ are empty,
$A^{(0)}$ is equal to $Q^{(0)}$ because of the first \kwfor loop, $Q^{(0)}$ contains only pairs $(K, 2)$ for every $K\in \cK$,
$\cL^{(0)}=\cK$, $n^{(0)}$ is not defined, $n^{(1)}=2$ and  $n^{(j)}\geq 2$ for $j\geq 1$.
In particular, properties~\ref{item:alg-a}-\ref{item:alg-value} are trivially satisfied for $l=0$.
Thus, fix $l > 0$ and assume the conditions are satisfied for $k<l$.

To see~\ref{item:alg-a} observe that in course of the $l$th iteration
of the \kwwhile loop we add an element to $A$  if and only if we add it to $Q$
and we never remove elements from $A$. This proves~\ref{item:alg-a}.

By the induction assumption $n^{(k)}\geq n^{(k-1)}$ for $k<l$. Hence,
in order to prove~\ref{item:alg-seq} we only need to show that  $n^{(l)} \geq n^{(l-1)}$.
The inequality is obvious when $n^{(l-1)}=2$. Thus, assume that $n^{(l-1)}>2$.
Then, $n^{(l-1)}$ and consequently also $n^{(l)}$ enter the queue $Q$ in pair with a respective knot
inside the \kwwhile loop. Hence,
there exist $k$ and $k'$ such that $l > k > k'$,
$ n^{(l)} = n^{(k)} + 1$ and $n^{(l-1)} = n^{(k')} + 1$.
From the induction assumption $n^{(k)} \ge n^{(k')}$, therefore
$ n^{(l)} = n^{(k)} + 1 \ge n^{(k')} + 1 = n^{(l-1)}$.

To prove~\ref{item:alg-uniq} let $l' < l$ be such that $(K^{(l)}, n^{(l)}) = (K^{(l')}, n^{(l')})$.
This obviously cannot happen if $n^{(l')}=2$. Thus,
for $l$ there exists a $k < l$ and for $l'$ there exists a $k' < l'$ such that
\[
  (K^{(l)}, n^{(l)}) \in Q^{(k)} \setminus Q^{(l)}
\]
and
\[
  (K^{(l')}, n^{(l')}) \in Q^{(k')} \setminus Q^{(l')}.
\]
Without loss of generality we may assume that $k' < k$. In the $k$th iteration line $(**)$ is executed with
$(K^{(l)}, n^{(l)}) \notin A^{(k)}$. However,
\[
(K^{(l)}, n^{(l)}) = (K^{(l')}, n^{(l')}) \in Q^{(k')} \subset A^{(k')} \subset A^{(k)},
\]
a contradiction.

Now we prove~\ref{item:alg-C}.
Assume to the contrary that there exists a
$
K \in \cK \setminus  \cK^{(l)}
$
such that $C^{(l)}[K]$ is not set.
Then, on the $l$th pass of the \kwwhile loop line $(*)$ is executed with $K^{(l)}=K$.
Indeed, if $K\not\in\cL^{(l-1)}$, then
$C^{(l-1)}[K]$ is set and the algorithm executes line $(*)$ by the induction assumption.
If $K\in\cL^{(l-1)}$, then $K=K^{(l)}$ and the algorithm executes the \kwelse branch, because  $C^{(l)}[K]$ is not set.
In particular, $(K,n)$ with $n=n^{(l)}$ is processed as one of the elements of $I^{(l)}[i]$ in the second $\kwforeach$ loop.
 If the algorithm executes line $(**)$, then $K\in\cK^{(l)}$, a contradiction.
 Otherwise $(K,n + 1) \in A^{(l-1)}$ and consequently $(K,n + 1) \in Q^{(l')}$ for some $l' < l$.
 Because of the assumption $(K, n+1)$ also is not in $Q^{(l)}$.
Hence,  $K=K^{(l'')}$ and $n+1 = n^{(l'')} > n^{(l)} = n$
for some $l''$ such that $l' < l'' < l$. This contradicts~\ref{item:alg-seq}.

Assume in turn that $C^{(l)}[K]$ is set and $K \in \cL^{(l)}$.
If $C^{(l-1)}[K]$ is set, then by the induction assumption $K\not\in \cL^{(l-1)}$.
Thus, the algorithm executes line $(**)$ in the $l$th iteration with $K_1=K$.
In particular, it also executes line $(*)$, a contradiction.
If $C^{(l-1)}[K]$ is not set, then $K^{(l)} = K$.
Let $n:=n^{(l)}$. Since $K\in\cL^{(l)}$ and the algorithm stops, there is a $k>l$ such that $K^{(k)}=K$.
Without loss of generality we may assume that $k$ is minimal with this property.
By~\ref{item:alg-seq} we have $n^{(k)}\geq n^{(l)}$ and by~\ref{item:alg-uniq} $n^{(k)} > n^{(l)}$.
It cannot be $n^{(k)} = n^{(l)}+1$, because on the $l$th pass of the \kwwhile loop the \kwelse branch is
not executed. Thus, $n^{(k)} > n^{(l)}+1$. It follows that there exists a $p<k$ such that $K^{(p)}=K$
and $n^{(k)} > n^{(p)}+1$. Therefore, $n^{(p)} > n^{(l)}$ and by~\ref{item:alg-seq} $p>l$.
Consequently, the number $k$, contrary to its choice, is not minimal.


Now we prove~\ref{item:alg-inj}.
It follows easily from the induction assumption that
if $C^{(l)}$ is not injective, then there exists a  $k < l$ such that
\[
  C^{(l)}[K^{(l)}] = C^{(l)}[K^{(k)}] = C^{(k)}[K^{(k)}].
\]
It means that
\[
 (n^{(l)}, i^{(l)}) =  (n^{(k)}, i^{(k)}).
\]
We set $C[K]$ only in the \kwthen branch of the first $\kwif$ statement.
Hence, $i^{(l)} \notin I^{(l-1)}$.
But in the $k$th iteration of the \kwwhile loop we set $C[K^{(k)}]$.
Therefore,  $i^{(k)} \in I^{(k)}$.
We do not delete elements from $I$.
Thus, since $k\leq l-1$, we have $i^{(l)} = i^{(k)} \in I^{(l-1)}$, a contradiction.

Again from the induction assumptions, if $D^{(l)}$ is not injective, then there exists a $k < l$ such that
\[
 K^{(l)} = D^{(l)}[i^{(l)}] = D^{(l)}[i^{(k)}] = K^{(k)},
\]
and $i^{(l)} \ne i^{(k)}$. Since $K^{(l)} = K^{(k)}$,
we have $i^{(l)}=I^n(K^{(l)}) = I^n(K^{(k)})=i^{(k)}$, a contradiction.

To prove~\ref{item:alg-value} let $i, K$ be such that $D^{(l)}[i]=K$.
Since $D$ is injective, we see that $K$ is unique.
If $D^{(l)}[i]=D^{(l-1)}[i]$ then \eqref{eq:alg-value} follows from  the induction assumption.
Otherwise, the algorithm sets $D$ in the $l$th iteration.
Thus, it also sets $C[K]$ to $(n, i)$.

To finish the proof first observe that assertion (a) is obvious.
Let $l^\star$ denote the final pass of the $\kwwhile$ loop.
Then, $Q^{(l^\star)}=\emptyset$. Thus, we get from~\ref{item:alg-C} that the keys of $C$ are all knots in $\cK$.
This proves assertion (b). Assertion (c) follows from~\ref{item:alg-value} and assertion (d) from~\ref{item:alg-inj}.
\qed

Given a family $\cK$ of prime knots and an integer $c$ we denote by $\cK_c$  the subfamily of these knots in $\cK$
whose minimal planar diagram requires $c$ crossings and we set $\cK_{\leq c}:=\bigcup_{i\leq c}\cK_i$.
Let $\HTW$ denote the family of $1~701~936$ prime knots tabulated by
Hoste, Thistlethwaite and Weeks in \cite{Hoste-etal-1998} and available from \cite{KnotAtlas}.
We have the following extension of Theorem~\ref{thm:cross-11}.

\begin{thm}
\label{thm:cross-14}
The invariant $I^n(K)$ distinguishes, up to mirror image,
between ambient isotopy classes of all prime knots in $\HTW_{\leq 14}$
and the classifying index of this family is $7$.
\end{thm}
\proof
   Apply Algorithm~\ref{alg:indexing} to the family $\HTW_{\leq 14}$.
   The algorithm stops, returning $7$ as the classifying index of the family.
   Hence the result follows from Theorem~\ref{alg:indexing}.
\qed

\begin{table}
\begin{tabular}{|c ||c |c |c |c |c |c |c |c |c |c |c |c |}
  \hline
  \hline
  $c$ & 3 & 4 & 5 & 6 & 7 & 8 & 9 & 10 & 11 & 12 & 13 & 14  \\
  \hline
  $N(c)$ & 2 & 2 & 3 & 3 & 3 & 3 & 5 & 5 & 6 & {\bf 6} & {\bf 7} & {\bf 7}   \\
  \hline
  \hline
\end{tabular}
\vspace*{3mm}
\caption{Classifying indexes $N(c)$ for the family $\HTW_{\leq c}$ with $c\leq 14$.}
\label{tab:classifying-indexes}
\end{table}

\begin{figure}
  \begin{center}
  \includegraphics[width=0.8\textwidth]{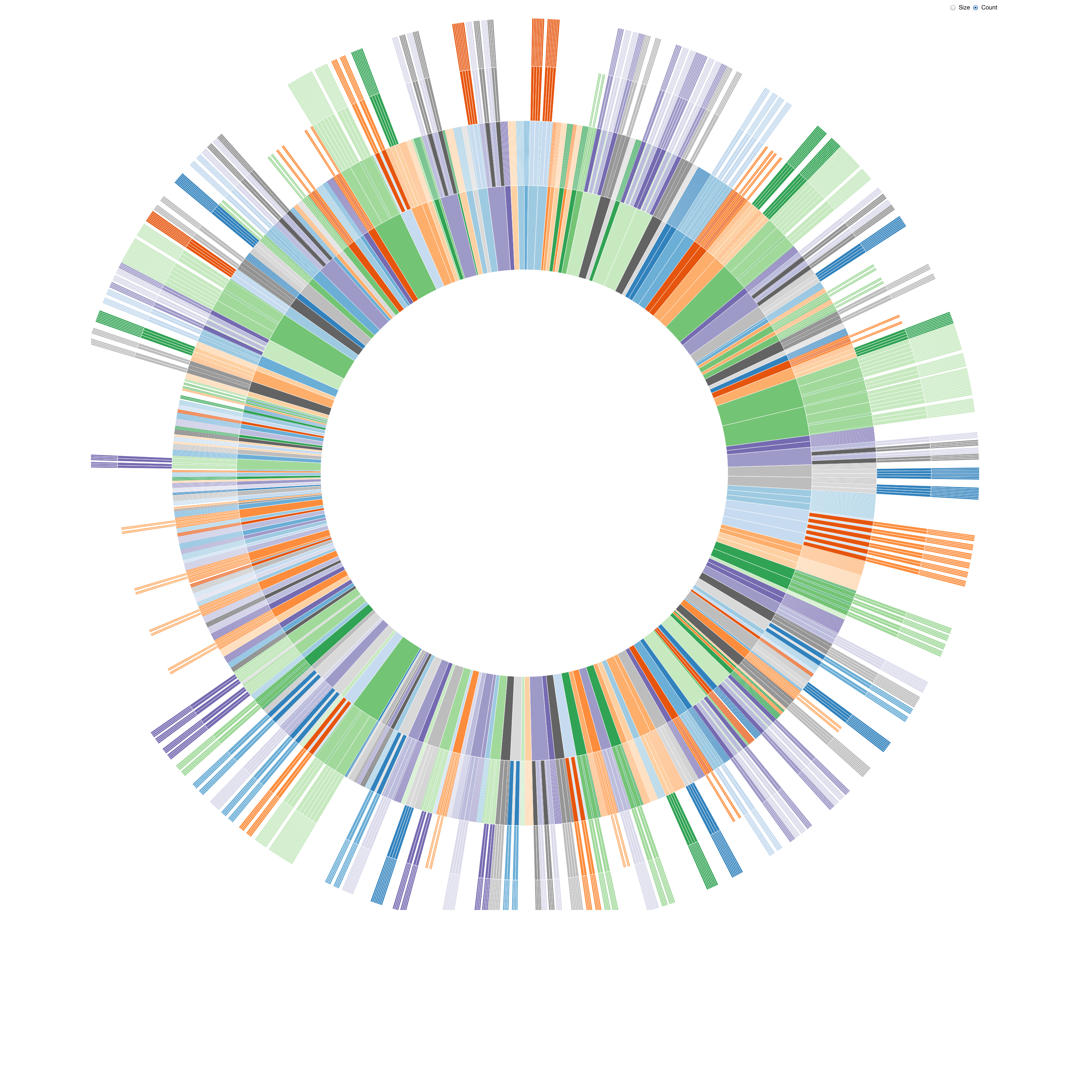}
  \end{center}
  \caption{The image shows knots separation in $\HTW_{\leq 14}$ by the invariant $I^n$. The middle
  white circle represents all knots in $\HTW_{\leq 10}$.
  The first ring  counting from inside presents arcs grouping knots with the same invariant with classifying index  $2$.
  The length of the coloured arc represents the group size. The following rings show groups for the classifying indexes $3$, $4$ and $5$. We can see how the groups split. Also,
  it is visible that the classifying index exceeds $3$ only for a small number of knots.}
  \label{fig:tree-10}
\end{figure}

The same way we can determine the classifying index for $\HTW_{\leq c}$
with $c\leq 14$. The results are presented in Table~\ref{tab:classifying-indexes}.
Figure~\ref{fig:tree-10} presents the grouping of classifying indexes of individual knots.
Tables~\ref{tab:classifying-index-12},~\ref{tab:classifying-index-13}
and~ \ref{tab:classifying-index-14} indicate that for most knots their classifying index
is relatively low.

\begin{table}
{\scriptsize
\begin{center}
\begin{tabular}{|l|lllll|}
\hline
 $c \backslash n$   & 2           & 3          & 4          & 5         & 6       \\
\hline
 3     & 1 (100\%)    & 1 (100\%)   & 0 (0\%)     & 0 (0\%)    & 0 (0\%)  \\
 4     & 1 (100\%)    & 1 (100\%)   & 1 (100\%)   & 0 (0\%)    & 0 (0\%)  \\
 5     & 2 (100\%)    & 2 (100\%)   & 2 (100\%)   & 2 (100\%)  & 0 (0\%)  \\
 6     & 3 (100\%)    & 3 (100\%)   & 3 (100\%)   & 2 (67\%)   & 0 (0\%)  \\
 7     & 7 (100\%)    & 7 (100\%)   & 6 (86\%)    & 3 (43\%)   & 0 (0\%)  \\
 8     & 21 (100\%)   & 21 (100\%)  & 14 (67\%)   & 10 (48\%)  & 1 (5\%)  \\
 9     & 49 (100\%)   & 49 (100\%)  & 30 (61\%)   & 19 (39\%)  & 1 (2\%)  \\
 10    & 165 (100\%)  & 165 (100\%) & 110 (67\%)  & 75 (45\%)  & 1 (1\%)  \\
 11    & 552 (100\%)  & 552 (100\%) & 355 (64\%)  & 225 (41\%) & 10 (2\%) \\
 12    & 2176 (100\%) & 2156 (99\%) & 1151 (53\%) & 727 (33\%) & 31 (1\%) \\
\hline
\end{tabular}
\end{center}
}
\vspace*{3mm}
\caption{Distribution of $I^n(K)$ computations for $\HTW_{\leq 12}$.
The entry in the $c$th row and $n$th column gives the number of knots $K\in\HTW_{c}$
for which it was necessary to compute the $I^n(K)$ in order to guarantee classification
in $\HTW_{\leq 12}$.
}
\label{tab:classifying-index-12}
\end{table}

\begin{table}
{\scriptsize
\begin{center}
\begin{tabular}{|l|llllll|}
\hline
 $n \backslash i$   & 2           & 3           & 4          & 5          & 6        & 7      \\
\hline
 3     & 1 (100\%)    & 1 (100\%)    & 0 (0\%)     & 0 (0\%)     & 0 (0\%)   & 0 (0\%) \\
 4     & 1 (100\%)    & 1 (100\%)    & 1 (100\%)   & 1 (100\%)   & 0 (0\%)   & 0 (0\%) \\
 5     & 2 (100\%)    & 2 (100\%)    & 2 (100\%)   & 2 (100\%)   & 0 (0\%)   & 0 (0\%) \\
 6     & 3 (100\%)    & 3 (100\%)    & 3 (100\%)   & 2 (67\%)    & 0 (0\%)   & 0 (0\%) \\
 7     & 7 (100\%)    & 7 (100\%)    & 7 (100\%)   & 3 (43\%)    & 0 (0\%)   & 0 (0\%) \\
 8     & 21 (100\%)   & 21 (100\%)   & 16 (76\%)   & 12 (57\%)   & 1 (5\%)   & 0 (0\%) \\
 9     & 49 (100\%)   & 49 (100\%)   & 40 (82\%)   & 24 (49\%)   & 1 (2\%)   & 0 (0\%) \\
 10    & 165 (100\%)  & 165 (100\%)  & 142 (86\%)  & 93 (56\%)   & 6 (4\%)   & 0 (0\%) \\
 11    & 552 (100\%)  & 552 (100\%)  & 474 (86\%)  & 301 (55\%)  & 15 (3\%)  & 0 (0\%) \\
 12    & 2176 (100\%) & 2176 (100\%) & 1692 (78\%) & 1004 (46\%) & 48 (2\%)  & 0 (0\%) \\
 13    & 9988 (100\%) & 9972 (100\%) & 7198 (72\%) & 4426 (44\%) & 293 (3\%) & 2 (0\%) \\
\hline
\end{tabular}
\end{center}
}
\vspace*{3mm}
\caption{Distribution of $I^n(K)$ computations for $\HTW_{\leq 13}$.
The entry in the $c$th row and $n$th column gives the number of knots $K\in\HTW_{c}$
for which it was necessary to compute the $I^n(K)$ in order to guarantee classification
in $\HTW_{\leq 13}$.
}
\label{tab:classifying-index-13}
\end{table}

\begin{table}
{\scriptsize
\begin{center}
\begin{tabular}{|l|llllll|}
\hline
 $n \backslash i$   & 2            & 3            & 4           & 5           & 6         & 7       \\
\hline
 3     & 1 (100\%)     & 1 (100\%)     & 1 (100\%)    & 1 (100\%)    & 1 (100\%)  & 0 (0\%)  \\
 4     & 1 (100\%)     & 1 (100\%)     & 1 (100\%)    & 1 (100\%)    & 0 (0\%)    & 0 (0\%)  \\
 5     & 2 (100\%)     & 2 (100\%)     & 2 (100\%)    & 2 (100\%)    & 0 (0\%)    & 0 (0\%)  \\
 6     & 3 (100\%)     & 3 (100\%)     & 3 (100\%)    & 2 (67\%)     & 1 (33\%)   & 0 (0\%)  \\
 7     & 7 (100\%)     & 7 (100\%)     & 7 (100\%)    & 5 (71\%)     & 2 (29\%)   & 0 (0\%)  \\
 8     & 21 (100\%)    & 21 (100\%)    & 21 (100\%)   & 13 (62\%)    & 2 (10\%)   & 0 (0\%)  \\
 9     & 49 (100\%)    & 49 (100\%)    & 46 (94\%)    & 26 (53\%)    & 3 (6\%)    & 1 (2\%)  \\
 10    & 165 (100\%)   & 165 (100\%)   & 158 (96\%)   & 105 (64\%)   & 9 (5\%)    & 1 (1\%)  \\
 11    & 552 (100\%)   & 552 (100\%)   & 523 (95\%)   & 329 (60\%)   & 26 (5\%)   & 0 (0\%)  \\
 12    & 2176 (100\%)  & 2176 (100\%)  & 2001 (92\%)  & 1253 (58\%)  & 94 (4\%)   & 0 (0\%)  \\
 13    & 9988 (100\%)  & 9988 (100\%)  & 8856 (89\%)  & 5494 (55\%)  & 451 (5\%)  & 3 (0\%)  \\
 14    & 46972 (100\%) & 46934 (100\%) & 38092 (81\%) & 23634 (50\%) & 2231 (5\%) & 21 (0\%) \\
\hline
\end{tabular}
\end{center}
}
\vspace*{3mm}
\caption{Distribution of $I^n(K)$ computations for $\HTW_{\leq 14}$.
The entry in the $c$th row and $n$th column gives the number of knots $K\in\HTW_{c}$
for which it was necessary to compute the $I^n(K)$ in order to guarantee classification
in $\HTW_{\leq 14}$.
}
\label{tab:classifying-index-14}
\end{table}

\subsection{Limitations of the method.}
Recall that during the computations of the $I^n$ invariant of a knot $K$
we first build a finite presentation of the fundamental group $G$ of $\RR^3\setminus K$
and construct the family $S^n(G)$ of the subgroups of index $n$ of $G$.

Since each subgroup of index $n$ of a finitely generated group $G$ corresponds to a homomorphism $G \to S_n$,
the number of such subgroups is bounded from above by $(n!)^g$ where $g$ is the number of generators in $G$.
We cannot bound the number from below. Thus, we can influence the computational time of
constructing $S^n(G)$ only via $g$.
Therefore, since $g$ is not greater than the number of critical cells of the constructed discrete vector fields,
it is important to minimize this number.
The problem of computing $S^n(G)$ is well studied in the literature. A good source of references is~\cite{Cosets}.

\begin{figure}
  \begin{center}
  \includegraphics[width=0.45\textwidth]{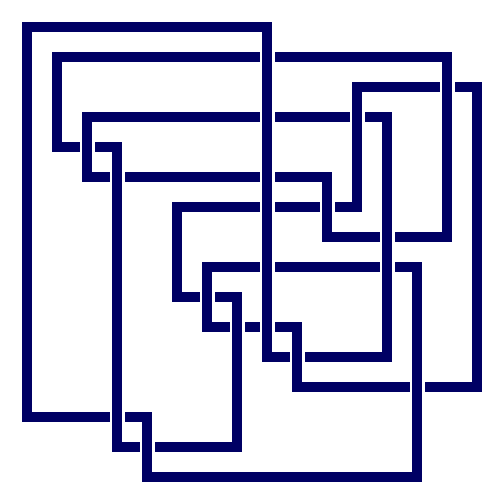}
  \end{center}
  \caption{Knot $\c{K}_{14,38437}$.}
  \label{fig:knot-14-38437}
\end{figure}

The problem of computing the optimal number of critical cells is NP-complete \cite{Joswig}.
In the construction of the discrete vector field we use greedy algorithms.
They usually behave very well for homology computation, where the difference
of one more or less critical cell has no significant consequence for the total computational complexity.
However, in the case of the $I^7$ invariant, one more critical cell may increase the total computation time
from minutes to days.
The most time-consuming knot we encountered is $\c{K}_{14,38437}$
presented in Figure~\ref{fig:knot-14-38437}.
For this knot the computations of $I^7(K)$ vary between a few minutes and $52$ days,
depending on the method used to construct the discrete vector field (coreductions, reductions)
and the order of cells in the data structure. The best method and order varies from knot to knot.
Since algorithm~\ref{alg:fundamentalGroup} proposed in this paper computes a group presentation
in a small fraction of time needed for the low-index subgroup computation, a reasonable strategy
to minimize the total computation time is to test a few methods for each knot in the search of the possibly small
number of critical cells.

The problem is additionally complicated by the Tietze transformations
used to simplify the group presentation. It turns out that the number of generators after transformations depends on  some
qualitative properties of the constructed discrete vector field. In Table~\ref{tab:gen-rel} we show the number of generators and relators before and after Tietze transformations in the case of four different strategies for the construction of the discrete vector field
for knot $\c{K}_{14,38437}$.
Only in one case we get $3$ generators and only in this case we are able to compute $I^7(K)$ in several minutes.
Our conclusion is that the number of generators does not determine the efficiency of the transformations.

Nevertheless, for the calculations it is important to perform shaving and  collapsibleSubset steps
even for knots given as arc presentation \cite{KnotAtlas}.
Without the geometric simplifications we cannot get less than $4$ generators for the knot in Figure~\ref{fig:knot-14-38437}. The steps are also important for applications, where knots are placed in a big cubical grid, e.g. 3D pictures of proteins.

\begin{table}
\begin{tabular}{||l|l|l||}
\hline
\hline
Method   & Before Tietze           & After Tietze     \\
\hline
 Reductions, cells order A   & g=4, r=4    &  g=4, r=3   \\
 Reductions, cells order B     & g=4, r=4    &  g=3, r=2   \\
 Coreductions, cells order A     & g=6, r=6    &  g=4, r=3   \\
 Coreductions, cells order B     & g=14, r=14    &  g=4, r=3   \\
\hline
\hline
\end{tabular}
\vspace*{3mm}
\caption{Number of generators (g) and relators (r) before and after Tietze transformations for
knot $\c{K}_{14,38437}$.
}
\label{tab:gen-rel}
\end{table}

\end{document}